\documentclass[12pt]{article}
\usepackage{graphicx}
\usepackage{mathptmx}
\usepackage{latexsym,amsmath,amssymb,amsfonts,amsthm}

\newcommand{\oM}{\overline{M}}
\newcommand{\ra}{\rightarrow}
\newcommand{\lnab}[1]{{\nabla}_{\!#1}}
\newcommand{\lnabb}[1]{\overline{\nabla}_{\!#1}}
\newcommand{\lnabo}[1]{{\nabla}^{\bot}_{\!\!#1}}

\newcommand{\ti}[1]{\mbox{\tiny $#1$}}
\newcommand{\scr}[1]{\mbox{\scriptsize $#1$}}
\newcommand{\ft}[1]{\mbox{\footnotesize $#1$}}
\newcommand{\sm}[1]{\mbox{\small $#1$}}
\newcommand{\la}[1]{\mbox{\large $#1$}}
\newcommand{\La}[1]{\mbox{\Large $#1$}}
\newcommand{\LA}[1]{\mbox{\LARGE $#1$}}
\newcommand{\Hu}[1]{\mbox{\Huge $#1$}}
\newcommand{\al}{\alpha}
\newcommand{\be}{\beta}
\newcommand{\h}{\mbox{\lie h}}                  
\newcommand{\g}{\mbox{\lie g}}
\newcommand{\R}{\mathbb{R}}          

\newfont{\lie}{eufm10 at 12pt}
\newfont{\field}{msbm10 at 11pt}

\newtheorem{theorem}{Theorem}[section]
\newtheorem{lemma}{Lemma}[section]
\newtheorem{corollary}{Corollary}[section]
\newtheorem{proposition}{Proposition}[section]

\begin{document}
\title{\normalsize  \bf   BERNSTEIN-HEINZ-CHERN RESULTS IN
CALIBRATED MANIFOLDS}
\author{Guanghan Li$^{(1)(2)}$ and Isabel M.C.\ Salavessa$^{(2)}$}
\date{}
\maketitle ~~~\\[-10mm]
{\small {\bf Abstract:} Given  a calibrated
Riemannian manifold $\oM$ with parallel calibration $\Omega$ of rank $m$
 and $M$  an orientable
m-submanifold with parallel mean curvature $H$, we prove that
if   $\cos\theta$ is bounded away from zero, where $\theta$ is
the $\Omega$-angle of $M$,
and if $M$  has zero Cheeger constant, then $M$ is  minimal.
 In the particular case $M$ is complete with
 $Ricc^M\geq 0$  we may replace
the boundedness condition on $\cos\theta$ by $\cos\theta\geq Cr^{-\beta}$,
when $r\ra +\infty$, 
where $0\leq \beta<1$ and $C>0$ are constants and $r$ 
is the distance function
to a point in $M$.
  Our proof 
is surprisingly simple and extends to a very large class of submanifolds in
calibrated manifolds, in a unified   way, the problem
started by Heinz and Chern of estimating
the mean curvature of graphic hypersurfaces in Euclidean
spaces. It is based on an estimation
of $\|H\|$ in terms of $\cos\theta$ and
an isoperimetric inequality. In a similar way,
we also give some conditions to conclude
$M$ is totally geodesic. We study  some particular cases.\\[2mm]
{\bf MSC(2000)}: 53C42; 53C38; 53C40; 58E35\\
{\bf Keywords:} calibrated geometry, parallel mean curvature, 
Heinz-inequality, Bernstein.
}

\markright{\sl\hfill G.\ Li and  I.\ Salavessa \hfill}

\section{Introduction}
\renewcommand{\thesection}{\arabic{section}}
\renewcommand{\theequation}{\thesection.\arabic{equation}}
\setcounter{equation}{0}
E.\ Heinz \cite{[H]} in 1955 introduced  the
problem of estimating the mean curvature of a surface of
$\R^3$ described by  a graph of a function
$f:\R^2\ra \R$.
He proved that if $f$ is defined on the disc $x^2+ y^2<r^2$ and
the mean curvature satisfies $\|H\|\geq c>0$, where $c$ is a constant,
then $r\leq \frac{1}{c}$. Thus, if $f$ is defined in all $\R^2$ and $\|H\|$
is constant, then $H=0$. Later, this result was extended 
 for the case of a map $f:\R^m \ra \R$ by Chern \cite{[Ch]} and
independently, by Flanders \cite{[F]}. This 
problem was generalized by the second author in her
Ph.D thesis (\cite{[S1]}, \cite{[S2]}) in 1987, for submanifolds
of  a Riemannian  product  $ \oM=M\times N$
of Riemannian manifolds $(M,g_1)$ and $(N,h)$, 
that can be
described as a  graph $\Gamma_{\!f}:=\{(p,f(p)):p\in M \}$ of a smooth map
$f:M\ra N$, that we recall as follows.
On any  oriented Riemannian manifold $(M,g)$ 
it is defined an isoperimetric
constant, the Cheeger constant
\begin{equation}
\h(M, g)=\inf_D\frac{A(\partial D,g)}{V(D,g)},
\end{equation}
where $D$ ranges over all open submanifolds of $M$ with compact
closure in $M$ and smooth boundary 
(see e.g.\ \cite{[Cha]} and section 4), and
$A(\partial D,g)$ and  $V(D,g)$ are respectively the area of 
$\partial D$ and the
volume of $D$, with respect to the metric $g$. We call such $D$ by
compact domain.
The Cheeger constant is zero, if, for example, $M$ is a closed manifold
(we abusively take the same definition for the closed case), or
if $M$ is a simple Riemannian manifold, that is, there exists
a diffeomorphism $\phi:(M,g)\ra (\R^{m}, <,>)$ onto $\R^m$
such that $\lambda^2 g \leq \phi^*<,>\leq \mu^2 g$ for some positive constants
$\lambda, \mu$. 
Another large class of Riemannian manifolds 
with zero Cheeger constant are the complete Riemannian
manifolds with non-negative Ricci tensor (see section 4).
Hence, zero Cheeger constant is a quite interesting condition.
\begin{theorem} (\cite{[S1]},\cite{[S2]}) 
If $f:(M,g_1)\ra (N,h)$ is a smooth map whose graph
 $\Gamma_{\!f}$ has parallel mean curvature $H$,
then for each   compact domain
$D\subset M$ we have the isoperimetric inequality\\[-3mm]
$$
\|H\|\leq {\frac{1}{m}}\frac{A(\partial D, g_1)}{V(D,g_1)}.$$
Thus $\|H\|\leq \frac{1}{m}\h(M,g_1)$.
In particular if $(M,g_1)$ has
zero Cheeger constant then $\Gamma_{\!f}$ is  a minimal
 submanifold of  $M\times N$.
\end{theorem}

We may also handle this  problem 
in the  context of calibrated manifolds.
A   calibration on a Riemannian manifold $\oM$ of dimension $m+n$
is a closed $m$-form $\Omega$ with comass one, that is, for each $p\in \oM$
and  any orthonormal system $X_i\in T_p\oM$,
$ |\Omega(X_1,\ldots ,X_m)|\leq 1 $ holds, 
and equality is achieved
at some system (see \cite{[HL]}).
If $F:M\ra \oM$ is an oriented immersed submanifold of dimension $m$, 
it is defined the  \em $\Omega$-angle \em of $M$,  
$\theta: M\ra [0,\pi]$, given by
$$\cos\theta=\Omega(X_1,\ldots, X_m),$$
where $X_i$ is a direct orthonormal frame of $T_pM$. 
We give to $M$ the induced metric $g=F^*\bar{g}$.
The submanifold is said to be  \em $\Omega$-calibrated  \em if 
$\cos\theta=1$. 
This is equivalent to $\Omega$ restricted to $M$ is 
the volume element of $M$.
Calibrated submanifolds are minimal, for they minimize
the volume of any domain $D$ among all variations
$F_t:D\ra \oM$, $t\in[0,1]$, of $F_0=F$ that fixes the boundary $\partial D$.
Let $dV_t$ be the volume element of $(D,g_t=F^*_t\bar{g})$.
Assuming $F_0$ is calibrated,
integration over $D$ of
$$
\cos\theta_1dV_1-dV_0=
F_1^*\Omega-F_0^*\Omega=d\tau
$$
 where
$\tau=\int_0^1 F_t^*(\Omega(\frac{\partial F}{\partial t}, \cdot))dt$ 
is a $(m-1)$-form that
satisfies $\tau_{|_{\partial D}}=0$,  gives
$$V_1(D)\geq \int_D\cos\theta_1dV_1= \int_D dV_0=V_0(D).$$
This inequality shows $F_0$ is minimal.
Furthermore, if $F_1$ also minimizes the volume
on the homotopy class of a calibrated submanifold, then $F_1$ is 
a calibrated submanifold as well. On the other hand, a stable minimal
submanifold $F$ may not be $\Omega$-calibrated. This is the
case $\oM$ has two different
$m$-calibrations and $F$ is calibrated only for one of them.
A pertinent question is to ask when is it true
that stable minimal submanifolds are in fact calibrated for some calibration.
This is true at least locally, for hypersurfaces 
in Euclidean spaces, or more generally for submanifolds
under certain  integrability conditions
 (see subsection 5.2).

The simplest examples of Riemannian manifolds with a calibration
are the Riemannian products $\oM=(M\times N,g_1\times h)$, 
with the volume calibration 
\begin{equation}
\Omega((X_1,Y_1), \ldots, (X_m,Y_m))= Vol_{(M,g_1)}(X_1,\ldots, X_m).
\end{equation}
If $M$ is a graph submanifold $\Gamma_{\!f}:M\ra M\times N$ 
then
$$\cos\theta= ({det(g_1 + f^*h)})^{-1/2}>0,$$ 
where the determinant is with respect to the metric $g_1$.
Reciprocally, a $m$-di\-mensio\-nal submanifold of $M\times N$
is locally  a graph if $\cos\theta>0$.
The graph is a calibrated submanifold if and only if $f$ is constant,
that is, the graph is a  slice. 
The condition $\cos\theta\geq \tau >0$, $\tau$ a constant,
 is equivalent to the boundedness
of $\|df\|^2$. The induced metric
on the graph $M$ is the graph  metric ${g}=g_1+f^*h$ on $M$
and so, under the above condition  the metrics $g$ and ${g}_1$ are
equivalent. In this case, 
$(M,{g})$ has zero Cheeger constant if and only if 
$(M,g_1)$ has so. 

In this paper
we will obtain the result in theorem 1.1 from a general result
for any  calibration $\Omega$, but with an extra condition
 on $\cos\theta$
at infinity. This means that this 
approach  for 
graphs is not so good has  the one in \cite{[S1]},\cite{[S2]}, 
although they are very much related to each other.
In both approaches we use a suitable vector field  $Z_1$, naturally
defined on all $M$ using the
calibration, but in theorem 1.1 we consider the
divergence of $Z_1$ with respect to the metric $g_1$ of $M$,
while in next theorem we consider the
divergence with respect to the induced metric $g$ of $M$.
On the other hand, we will provided a unified  way to obtain
a Heinz-Chern result for  submanifolds with
parallel mean curvature in a very large class of ambient spaces, the class of
calibrated manifolds.

Examples of calibrated manifolds are the K\"{a}hler manifolds with
the K\"{a}hler calibration, the
Riemanniam manifolds
with special holonomy, namely, the Calabi-Yau manifolds
with the special Lagrangian calibration, the
qua\-ter\-nio\-nic-K\"{a}hler manifolds with the quaternionic calibration, 
the hyper-K\"{a}hler  manifolds (with many calibrations), 
$G_2$ manifolds 
with the associative and co-associative calibration, and $Spin(7)$
manifolds with the Cayley calibration (see  
\cite{[Jo]}).
These special spaces are Einstein  manifolds, and except
the quaternionic-K\"{a}hler case, they are all Ricci flat.
If 
 $n=1$, a parallel $\Omega$ defines a non-zero
 global parallel vector field
$(*\Omega)^{\sharp}$ on $\oM$
and so, if $\oM$ is simply connected
then it splits as a Riemannian product
$\oM=M\times N^1$, where $N^1$ is one dimensional, and
$\Omega$ is the volume element of $M$.  More generally, 
for $n=1$, a divergence free vector field $\bar{X}$ on $\oM$,
 defines a closed $m$-form 
$\Omega=*\bar{X}^{\flat}$, where $*$ is the star operator
on $\oM$. This form is a calibration if $\|\bar X\|=1$. This is the case
of a Riemannian manifold $\oM$ with a codimension-one 
transversally oriented foliation by minimal
hypersurfaces, for, in this case the unit normal $\bar{X}$ to the leaves 
defines a divergence free vector field  of $\oM$ that calibrates the
leaves. For foliations of any codimension  see section 5.
\\[-1mm]

In what follows, $(\oM, \bar{g},\Omega)$ denotes a 
calibrated  $(m+n)$-dimensional manifold with a
calibration $\Omega$ of rank $m\geq 2$,
and $F:M\ra \oM$ is an immersed oriented submanifold
of dimension $m$, induced metric $g$, volume element
$dV$, normal bundle $NM$,
mean curvature $H$ and
$\Omega$-angle $\theta$. 
We consider the following 
morphisms $\Phi:TM\to NM$, $\Psi:\wedge^2TM\to\wedge^2NM$, such that 
  for $X,Y\in T_pM$, $U,V\in NM_p$,
\begin{equation}
\begin{array}{lcl}
\bar{g}(\Phi(X),U) &=&\Omega(U,*X),\\
\langle \Psi(X\wedge Y),U\wedge W\rangle &=&
\Omega(U\wedge W, *X\wedge Y), 
\end{array}
\end{equation}
where $*:TM\to\wedge^{m-1}TM$ and $*:\wedge^2TM\to \wedge^{m-2}TM$ are 
the star operators and $\langle,\rangle$ is the usual inner product
in $\wedge^2NM$. For $m=2$, set $\langle \Psi(dV),U\wedge W\rangle
=\Omega(U\wedge W)$, where $dV$ is the volume element of $M$.
Our main results are:
\begin{theorem}[The integral $\Omega$-isoperimetric inequality]
On a compact domain $D$ of $M$,
with boundary $\partial D$ with volume element $dA$,
 the following inequality holds
$$ \left|\int_D \LA{(}-m\cos\theta \|H\|^2+\langle \nabla{}^{\bot}H,
\Phi\rangle \LA{)} dV
+\int_D\bar{\nabla}_{H} \Omega~\right|\leq 
\int_{\partial D}\sin\theta\, \|H\|dA, 
$$
where $\langle,\rangle$ is the Hilbert-Schmidt inner product in
$TM^*\otimes NM$. 
\end{theorem}
From now on we assume
$\Omega$ is parallel. Theorem 1.2 leads to:
\begin{theorem} 
If $F:M\ra \oM$ is immersed   with parallel mean curvature
and 
$\cos\theta>0$,  on a compact domain $D$ of $M$,
the following  isoperimetric inequality holds:
$$
\|H\|~\leq ~\frac{1}{m}
\,\,\left( \frac{\sup_{\partial D}\sin\theta}
{ \inf_D \cos\theta }\right)
\,\, \frac{A(\partial D,g)}{V(D,g)}.$$
 In particular:\\[1mm]
(1) If 
 $\cos\theta\geq \tau>0$ where
$\tau$ is a constant, then 
$\| H\|\leq \frac{{1}}{m}\, \, \frac{\sqrt{1-\tau^2}}{\tau}\, \, \h(M,g).$
In this case, if $M$ has zero Cheeger constant, then $M$ is a minimal 
submanifold.\\[1mm]
(2) If $\cos\theta=1$ on $\partial D$ for some domain $D$ then
$F$ is a minimal immersion.
\end{theorem}
\begin{corollary} If $M$ is closed with parallel mean curvature
and $\cos\theta>0$, then $M$ is minimal.
\end{corollary}
\begin{corollary} If $M$ is closed, 
 $1> |\cos\theta|$ constant, and  $\|\Phi(X)\|\leq \mu \sin\theta\,\|X\|$,
where $0< \mu\leq 1$ is a constant, and $\|H\|$ not identically zero,
then
$$|{\cot}\, \theta |
\leq \frac{\mu}{\sqrt{m}}\frac{\int_M\|\nabla^{\bot}H\|dV}{\int_M \|H\|^2dV}.$$
Equality holds iff  $\Phi$
is a homothetic
morphism with coefficient of conformality $\mu^2\sin^2\theta$ 
on the orthogonal complement  of the distribution defined by the
kernel of $\nabla^{\bot}H$, and 
$\nabla^{\bot}_{X}H= \Phi(\psi(X))$
where $\psi:TM\to TM$ is a linear morphism.
\end{corollary}
\noindent
We will see that $\|\Phi(X)\|\leq \sin\theta \, \|X\|$ always hold.
The conformality condition on $\Phi$ is not an uncommon condition.
In a 8-dimensional quaternionic K\"{a}hler
manifold,
almost complex 4-submanifolds define a morphism
 $\Phi$ with  coefficient of conformality
$(1-\cos\theta)(\cos\theta-\frac{1}{3})$.
Four dimensional submanifolds with equal K\"{a}hler angles of a 
K\"{a}hler manifold of complex dimension 4, 
define  $\Phi$ with coefficient of conformality $
(1-\cos\theta)\cos\theta$ (see section 5).\\[-3mm] 

 We can slightly improve theorem 1.3 in case $Ricci^M\geq 0$
and $M$ is complete. In this case, if we fix $p\in M$, there is a constant
$C_1>0$, such that (see section 4)
\begin{equation}
\h(M)\leq \h(B_r(p))\leq \frac{C_1}{r}~~~~~~~\mbox{for~all~} r\in (0,+\infty).
\end{equation}
\begin{theorem} 
If $F:M\ra \oM$ is a complete  immersed oriented $m$-dimensional
 submanifold with parallel mean curvature, and 
$Ricci^M\geq 0$ 
and the $\Omega$-angle satisfies
$\cos\theta\geq  Cr^{-\beta}> 0$ when $r\rightarrow +\infty$,
where $0\leq \beta<1$ and $C>0$ are constants, and $r$ is the distance
function in $M$ to a point $p\in M$, 
then $F$ is a minimal submanifold.
\end{theorem}
An application of theorem 1.1 is the following:
\begin{proposition} If $(M,g_1)$ is a complete Riemannian manifold
with $Ricci^{(M,g_1)}\geq 0$, then any graphic submanifold
with parallel mean curvature
 $F=\Gamma_{\!f}:M\ra (M\times N, g_1\times h)$, where
$f:(M,g_1)\ra (N,h)$ is a smooth map, 
  is a minimal submanifold.
\end{proposition}

 It is fundamental 
some nonnegativeness on the curvature tensor of $M$ 
 to obtain such Heinz-Chern results. 
If $\oM=\mathbb{H}^m\times \mathbb{R}$ where $\mathbb{H}^m$ is 
 is the $m$-hyperbolic space, there are  examples
of entire graphic hypersurfaces, and so complete,
 with non-zero constant mean curvature $c$ and with
$\cos\theta$ bounded away from zero,  
as can be shown by the following proposition. 
Note that $\h(\mathbb{H}^m)=m-1$.
The function
$r(x)=\ln \left({(1+|x|)}/{(1-|x|)}\right)$
is the distance function
in $\mathbb{H}^m$ to $0$, for the Poincar\'{e} model,
and $\nu=(-\nabla f,1)/\sqrt{1+\|\nabla f\|^2}$ is a unit normal
to $\Gamma_{\!f}$:

  \begin{proposition} \cite{[S1],[S2],[S3]} 
For each $ |c| \leq m-1$,  $f_c:\mathbb{H}^m\ra \R$ defined by:
\[f_c(x)=\Hu{\int}_0^{r(x)}
\frac{\frac{c}{(\sinh r)^{m-1}}\int_0^r(\sinh t)^{m-1}dt}
{\sqrt{1-\La{(}\frac{c}{(\sinh r)^{m-1}}\int_0^r(\sinh t)^{m-1}dt\La{)}^2}}
dr, \]
 is smooth on all $\mathbb{H}^m$, and for
each $d\in \R$, $\Gamma_{f_c +d}\subset
\mathbb{H}^m\times \R$ has constant mean curvature given
 by $\bar{g}(H,\nu)
=\frac{c}{m}$, and $\cos\theta > \sqrt{\scriptsize{(m-1-|c|)/(m-1)}}$.
Furthermore, 
$\{\Gamma_{(f_c)+d}(x):~ x\in \mathbb{H}^m, d\in \mathbb{R}\}$ 
and $\{\Gamma_{(f_c)+d}(x):~ x\in \mathbb{H}^m, c\in [1-m, m-1]\}$
define (partial) foliations of $\mathbb{H}^m\times \mathbb{R}$
 by hypersurfaces with
the same constant mean curvature $c$, and with constant mean curvature
parameterized by the leaf, respectively.
\end{proposition}
A  related  classical problem  is
the  Bernstein-type problem, that determines when a minimal submanifold
must be totally geodesic.
In 1927, Bernstein \cite{[Bern]} proved that any minimal surface of $\R^3$
defined by the graph of an entire  map $f:\R^2\to\R$ is a linear plane.
This result was generalized to $\R^{m+1}$ for $m\leq 7$
by de Giorgi \cite{[deGio]} ($m=3$), Almgren \cite{[Alm]}
($m=4$), and Simons \cite{[Simons]}($m\leq 7$), 
and to
higer dimensions and codimensions under various growth conditions
by many others, as for example  Hildebrandt, Jost,
and Widmann in \cite{hjw}, 
Ecker and Huisken \cite{eh}, Wang \cite{wa2}, and more recently 
some attention is given to
Bernstein theorems in curved
Riemannian product or warped product spaces
by Al\'{\i}as, Dajczer and Ripoll \cite{[ADR]}.
In higer dimension, and mainly in higer codimension, Bernstein-type results
tend to be more difficult and complicated to formulate. 
Some Bernstein results have been  obtained 
for stable minimal hypersurfaces by do Carmo and Peng \cite{[doCarmoPeng]}, 
Miranda \cite{[Mira]}, Fischer-Colbrie, Schoen, 
Simon and Yau \cite{fcs,[SSY]}, and
 for leaves of transversely oriented codimension one
foliations of Riemannian manifolds by Barbosa, Kenmotsu, Oshikiri, Bessa
and Montenegro \cite{[BKO],[BBM]},  
where Chern-Heinz inequalities are derived, as well the stability
of the leaves.\\[-3mm]

In this paper we obtain some Bernstein-type results using the same
philosophy of the Chern-Heinz inequalities, applied to  submanifolds
immersed in calibrated manifolds, and under certain conditions,
 allowing us to obtain this type of
results in any codimension. They are derived from the expression of
 $\Delta\cos\theta$.
This Laplacian involves the covariant derivative of the mean curvature,
a quadratic term on  the second fundamental form $B$ and a curvature term of
$\oM$  that we have to analyse. We should have in mind that if
 $F$ is totally geodesic then
$\theta$ is constant. 
Let $B_{\Phi}$ a 1-form on $M$ and
$\wedge^2 B:\wedge^2TM \to \wedge^2NM$ given by
\begin{equation}\begin{array}{rcl}
B_{\Phi}(X) &=& \sm{\sum}_i\, \bar{g}(B(X_i,X),\Phi(X_i))\\
\wedge^2 B(X \wedge Y) &=& \sm{\sum}_i \, B(X_i,X)\wedge B(X_i,Y),
\end{array}
\end{equation}
where $X_i$ is any  orthonormal basis of $T_pM$.
We   consider the following quadratic forms
 defined for any $m$-calibration $\Omega$ and $F:M^m\to \oM$ satisfying
$\cos\theta>0$, and  applied to
tensors $B'\in \bigodot^2TM^*\otimes NM$
\begin{equation}
 \begin{array}{l}
Q_{\Omega}(B')=\tilde{Q}_{\Omega}(B')
+\frac{1}{\cos^2\theta}\|B'_{\Phi}\|^2\\[2mm]
\tilde{Q}_{\Omega}(B') = \|B'\|^2-\frac{2}{\cos\theta}\langle
\Psi,\wedge^2 B'\rangle. \end{array}
\end{equation}
 $Q_{\Omega}$ (or $\tilde{Q}_{\Omega}$)  is said to be
\em  $\delta$-positive \em at $B'$ if 
$Q_{\Omega}(B')\geq \delta \|B'\|^2
$, where $\delta>0$. This is the case, with $\delta=1$,
 when $\langle \Psi,\wedge^2 B'\rangle
\leq 0$, as it is the case $n=1$. 

The quadratic form $Q_ {\Omega}$ was defined in \cite{wa2} for the case
$\oM=M^m\times N^n$ and $\Omega$  the volume element of $(M,g_1)$
and $F=\Gamma_f$ with
$f:(M,g_1)\to (N, h)$ a smooth map.
For   $n\geq 2$, one has  to require 
$\lambda_i\lambda_j\leq (1-\delta)$, for $i\ne j$ 
and  some constant $0<\delta\leq 1$,
where $\lambda_1^2\geq \ldots\geq \lambda_m^2\geq 0$ 
are the eigenvalues of $f^*h$,  to have
$Q_{\Omega}$ $\delta$-positive at any
 $B'$. This  condition
 gives bounds to the components of the calibration $\Omega$ in a
convenient basis $X_i$ of $TM$ and $U_{\alpha}$ of $NM$, namely on the
components
$\langle\Psi(X_i\wedge X_j), U_{\alpha}\wedge U_{\beta}\rangle= 
\cos\theta\delta_{i\alpha}\delta_{j\beta}\lambda_i\lambda_j$
 (see section 5).
In general, this condition on $Q_{\Omega}$
or  on $\tilde{Q}_{\Omega}$, can be
 a quite restrictive condition for the higer codimension case, and it 
holds for calibrated submanifolds only if these
are  necessarily totally geodesic (see proposition
1.3). But it holds for certain kind of submanifolds, as
 for example, if $F$ is sufficiently close to a  totally umbilical
submanifold (i.e. satisfy $B=Hg$).
For  the K\"{a}hler calibration we will find in proposition 5.3
that $\delta$-positiveness is quite unlike to hold
on minimal 4-submanifolds with equal K\"{a}hler angles unless $\oM$
is a Calabi-Yau 4-fold.
In lemma 5.5 for almost complex submanifolds
in quaternionic K\"{a}hler manifolds, we give a natural condition
that ensures $\delta$-positiveness of $Q_{\Omega}(B)$.\\ 

For calibrated submanifolds we have:
\begin{proposition} If $F:M\to \oM$ is $\Omega$-calibrated then 
$Q_{\Omega}(B)=\tilde{Q}_{\Omega}(B)=0$.
Thus, if  ${Q}_{\Omega}$  
is $\delta$-positive, then
$F$ is totally geodesic. This is always the case $n=1$.
\end{proposition}
\noindent
\begin{theorem} Assume $M$ has parallel mean curvature,
 $\cos\theta>0$, and 
$Q_{\Omega}(B)\geq \delta \|B\|^2$ for some constant $\delta>0$ and
\begin{equation} 
\sm{\sum}_{ij}\bar{R}(X_i,X_j,X_i,\Phi(X_j))\geq 0.
\end{equation} 
Then the following inequalities hold for any compact domain $D$:
\begin{equation}\begin{array}{c}
\|\nabla \log\cos\theta\| \leq  \sqrt{m} \tan\theta\, \|B\|,
\\[1mm]
 \int_{D}\|B\|^2dV  \leq  \frac{\sqrt{m}}{\delta}\int_{\partial D}
\tan\theta\, \|B\| dA.
\end{array}\end{equation}
Moreover:\\[1mm]
(A) if ~$\tan\theta\, \|B\|$ is integrable on $M$ and 
$M$ is complete, then $F$ is 
totally geodesic;\\[1mm]
(B) if ~$\inf_M\cos\theta=\tau>0$,
and $\|B\|$ is not identically zero, then 
$$\frac{\inf_M\|B\|^2}{\sup_{M}\|B\|}\leq
\inf_D\left(\frac{-\!\!\!\!\!\!\!\int_D \|B\|^2dV}
{-\!\!\!\!\!\!\!\int_{\partial D}\|B\| dA}\right)
\leq \frac{\sqrt{m}}{\delta}\frac{\sqrt{1-\tau^2}}{\tau}\, \h(M).$$
In particular, if $\|B\|$ is constant, then $\|B\|\leq \frac{\sqrt{m}}
{\delta}\frac{\sqrt{1-\tau^2}}{\tau}\h(M)$. In this case (since $\|B\|\neq 0$),
$\h(M)\neq 0$.\\[1mm]
(C) If the sectional curvatures of $\oM$ are bounded from below,
$M$ is complete, 
$\|B\|$ is  bounded,  then either $\inf_M\cos\theta =0$ or $\inf_M\|B\|=0$.
\end{theorem}
\noindent
$\delta$-positiveness of $Q_{\Omega}(B)$ and (1.7) imply
$\Delta\log\cos\theta\leq -\delta\|B\|^2$, 
and (A)-(C) are  consequences of this.
 Parallel submanifolds (i.e.\ $\nabla B=0$) have parallel mean curvature,
$\|B\|$ is constant, and equality to zero at (1.7). 
We will see in section 5 that
(1.7)  holds for example for  
 almost complex submanifolds
in  quaternionic space forms with nonnegative
scalar curvature.
If $n=1$, (1.7)
is equivalent to $\sum_jRicci^{\oM}(X_j,\Phi(X_j))\geq 0$, that  holds
for $\oM$ an Einstein manifold.\\

For $M$ noncompact surface, using a criteria for parabolicity we  prove  next
theorem:
\begin{theorem} Assume $F:M^2\to \oM$ is a  minimal complete
 immersed surface  with 
$\cos\theta>0$, and  $\bar{K}\circ F\geq 0$ away from a compact set of $M$,
where $\bar{K}$
denotes the sectional curvatures of $\oM$.
 If (1) or (2) below holds:\\[1mm]
(1)  $\oM$ is a space form of dimension 3;
\\
(2) 
$Q_ {\Omega}(B)\geq \delta \|B\|^2$, $\cos\theta\geq \tau$, with
$ \tau, \delta>0$ constants  and (1.7) holds;\\[1mm]
then $M$ is a totally geodesic submanifold.
\end{theorem}
\noindent
The previous theorem applied to graphs gives (simpler proofs) of the classical
Bernstein results:
\begin{corollary}
(1) \cite{[Bern],[Ch2]}
If a smooth entire function $f:\mathbb{R}^2\to \mathbb{R}$ defines 
a minimal graph in
$\mathbb{R}^3$, then $f$ is linear.\\[1mm]
(2) \cite{wa2} 
If  a smooth entire function $f:\mathbb{R}^2\to \mathbb{R}^{n}$ defines 
a minimal graph in $\mathbb{R}^{n+2}$, with $\|\wedge^2df\|
=|\lambda_1\lambda_2|
\leq 1-\delta$, 
$1\geq \delta >0$  constant, and  $ \|df\|$ bounded 
(or equivalently, $\cos\theta\geq \tau> 0 $),
then  $f$ is a linear map.
\end{corollary}
\noindent
An application of theorem 1.5 gives the following
 Bernstein-type results for graphic submanifolds:
\begin{corollary}
(3) Let $f:(M^m,g_1)\to (N^n,h)$ defining a minimal graph on
$\oM=M\times N$ with $ \int_M\tan\theta\, \|B\|dV<+\infty$. 
We assume $(M,g_1)$ is complete
with sectional curvature $K_1$ and Ricci tensor $Ricci_1$
and $(N^n,h)$ has sectional curvature $K_N$ satisfying:
 (a) for $n=1$,
 $Ricci_1\geq 0$; 
 (b) for $n\geq 2$,
 $f^*h<(1-\delta)g_1$, $1\geq \delta>0$, 
 and  at each $p\in M$ and two-planes $P$ of $T_pM$,
and $P'$ of $T_{f(p)}N$,
 either $K_1(P)\geq K_N(P')^+$,
or  $Ricci_{1}(p)\geq 0$ and  
$K_N(P')\leq -K_{1}(P)$. Then
 $f$ is totally geodesic.
Furthermore if at some point $K_{1}(p)>0$ (or $Ricci_{1}>0$)
then $f$ is constant, and $\Gamma_f$ is a slice.
\end{corollary}
\noindent
For immersed  hypersurfaces in the Euclidean space we have:
\begin{corollary} Assume $M$ is a complete minimal immersed  hypersurface of
$\mathbb{R}^{m+1}$,
such that  for some parallel calibration
$\Omega$ of ~$\mathbb{R}^{m+1}$,  we have $\cos\theta>0$ on all $M$
(and so, $M$ is locally a graph).
If   $\int_M \sin\theta\|B\|dV<+\infty$, 
then $M$ is a linear hyperplane. 
\end{corollary}
In section 2 we derive the fundamental properties of
the morphism $\Phi$ and  prove theorem 1.2 and corollary 1.2.
In section 3 we discuss when $F$ is totally geodesic,
describe the formula of the Laplacian of $\cos\theta$ and 
give the proof of proposition 1.3,  theorems 1.5 and 1.6, and corollaries 1.3,
 1.4 and 1.5.
In section 4 we obtain some properties of the Cheeger constant,
 and prove theorem 1.4 and 
proposition 1.1.
In section 5
we specify to some examples of submanifolds in calibrated Riemannian manifolds,
namely in a foliated space, in 
K\"{a}hler and  quaternionic-K\"{a}hler
manifolds. 
\section{The morphism $\Phi$}
\renewcommand{\thesection}{\arabic{section}}
\renewcommand{\theequation}{\thesection.\arabic{equation}}
\setcounter{equation}{0}
Given $\Omega$ a $m$-calibration on $\oM$, we
consider the  $T\oM$-valued $(m-1)$-form
$\Omega^{\sharp}:\wedge^{m-1}T\oM^*\ra T\oM$,
$\bar{g}(\Omega^{\sharp}(X_2\,\ldots,X_{m}), X_1)=\Omega(X_1,
\ldots, X_m)$,
where $X_i \in T_p\oM$.
Let $F:M\ra \oM$ be an immersed submanifold of dimension $m$
with normal bundle $NM$, and $\Omega$-angle $\theta$.
We denote by $\lnabb{}$, $\lnab{}$ and $\lnabo{}$ the respective
covariant derivatives of $\oM$, $M$ and $N
M$,
and by $B(X,Y)=\lnab{X}dF(Y)$ the second fundamental form of $F$,
defined by the following equations for $X,Y$ vector fields on $M$ and
$U$ section of $NM$
$$ \lnab{X}Y=(\lnabb{X}Y)^{\top},~~~~~~(\lnabb{X}Y)^{\bot}=B(X,Y),~~~~~~
\lnabo{X}U= (\lnabb{X}U)^{\bot},$$
where $(\cdot)^{\top}$ and $(\cdot)^{\bot}$ are the orthogonal projections
into $TM$ and $NM$ respectively.
The mean curvature is $H=\frac{1}{m}trace\, B$.
We consider the
 morphism  $\Phi=\Phi_{\Omega}:TM\ra NM$ defined in (1.3),
$\Phi(X)=(\Omega^{\sharp}(*X))^{\bot}$.
Recall the covariant derivative and the co-differential
of $\Phi$ are given by
$$\lnab{X} \Phi(Y)=\lnabo{X}(\Phi(Y))-\Phi(\lnab{X}Y), ~~
~~~~\delta \Phi=-\sum_i \lnab{X_i}\Phi(X_i).$$
\begin{lemma} If $X\in T_pM$ and $U\in NM_p$ are units, then
$|\bar{g}(\Phi(X),U)|\leq \sin\theta$.
\end{lemma}
\noindent
\em Proof. \em Let  $X_i$ be a direct o.n.\ basis 
of $T_pM$ with $X_1=X$.
Consider the function
$$\phi(t)=\Omega(\sm{\frac{\epsilon X_1+tU}{\sqrt{1+t^2}}},X_2,\ldots,X_m)=
\sm{\frac{1}{\sqrt{1+t^2}}}|\cos\theta| +\sm{\frac{t}{\sqrt{1+t^2}}}
\bar{g}(\Phi(X),U)$$
where $\epsilon=\pm 1$ s.t.\ $\epsilon \cos\theta=| \cos\theta|$.
Since $\Omega$ is a calibration, $\phi(t)\leq 1$ for any $t$.
We may assume $\cos\theta\neq 0$. 
At $t=\frac{\bar{g}(\Phi(X),U)}
{|\cos\theta|}$,  
$\phi(t)=\sqrt{\cos^2\theta+ \bar{g}(\Phi(X),U)^2}\leq 1$. \qed
\begin{lemma} For $X_i$ and $U_{\alpha}$   d.o.n. basis 
of $T_pM$ and $NM_p$ respectively, we have
\begin{eqnarray*}
\delta \Phi &=& m\, \cos\theta \, H 
-\sm{\sum}_{\alpha}(\overline{\nabla}_{U_{\alpha}}\Omega)
(X_1,\ldots,X_m)U_{\alpha}\\ \langle\nabla\Phi,B\rangle &=&
-\cos\theta\, \tilde{Q}_{\Omega}(B)+
\sm{\sum}_{jk}(\overline{\nabla}_{X_j}\Omega)(X_i,\ldots, B(X_j,X_k)_{\scr{(k)}},
\ldots, X_m)
\end{eqnarray*}
\end{lemma}
\noindent
\em Proof. \em Let $U$ a section of $NM$ and 
 $Z\in T_{p}M$.  We may assume
 $\lnabo{}U(p)=\lnab{}X_i(p)=0$. At $p$,
$g ((\lnabb{Z}U)^{\top}, X) =-\bar{g}( B(Z,X),U)$, 
$\lnabb{Z}X_i(p)=
(\lnabb{Z}X_i(p))^{\bot}=B(Z,X_i),$ and
\begin{eqnarray*}
\bar{g}( \lnab{Z}\Phi(X_1),U)
&=& Z\cdot \bar{g}( \Phi(X_1),U)
 =d(\Omega(U, X_2,\ldots, X_m))(Z)\\
&=& \overline{\nabla}_{Z}\Omega(U,X_2,\ldots,X_m)+
\Omega((\lnabb{Z}U), X_2, \ldots, X_m)\\
&&+\sm{\sum}_{i\geq 2}
\Omega(U, X_2,\ldots,\lnabb{Z}X_i, \ldots, X_m)\\
&=&\overline{\nabla}_{Z}\Omega(U,X_2,\ldots,X_m) +
\Omega((\lnabb{Z}U)^{\top}, X_2,\ldots, X_m)\\
&&+\sm{\sum}_{i\geq 2}
\bar{g}\la{(} \Omega^{\sharp}( X_2,\ldots, B(Z,X_i), \ldots, X_m), 
U\la{)}.
\end{eqnarray*}
That is
\begin{eqnarray*}
\lnab{Z}\Phi(X_1) &=& 
\sm{\sum}_{\alpha}\overline{\nabla}_{Z}\Omega(U_{\alpha},X_2,\ldots,X_m)
U_{\alpha} -\cos\theta \, B(Z,X_1)\\
&& +\sm{\sum}_{i\geq 2} (\Omega^{\sharp}( X_2,\ldots,
B(Z,X_i),\ldots, X_m))^{\bot}.
\end{eqnarray*}
 We have
$*X_k=(-1)^{k-1}X_1\wedge \ldots \wedge \hat{X}_k \wedge \ldots
\wedge X_m$. Hence,
\begin{eqnarray*}
\lnab{Z}\Phi(X_k) &=&  \sm{\sum}_{\alpha}(-1)^{k+1}
\overline{\nabla}_{Z}\Omega(U_{\alpha},X_1,\ldots,\hat{X}_k,\ldots
X_m)U_{\alpha} -\cos\theta \,B(Z,X_k)\\
&&+
\sm{\sum}_{1\leq i<k}(-1)^{k+1} (\Omega^{\sharp}( X_1, \ldots, B(Z,X_i), 
\ldots, \hat{X}_k,\ldots,  X_m))^{\bot}\\
&&+\sm{\sum}_{k<i} (-1)^{k+1}(\Omega^{\sharp}( X_1,\ldots,  \hat{X}_k,\ldots,
B(Z,X_i), \ldots, X_m))^{\bot}.
\end{eqnarray*}
Therefore,
\begin{eqnarray*}
\sm{\sum}_k\lnab{X_k}\Phi(X_k) &=&
\sm{\sum}_{\alpha} -d\Omega(U_{\alpha},X_1,\ldots,X_m)U_{\alpha}\\
&&+\sm{\sum}_{\alpha}
(\overline{\nabla}_{U_{\alpha}}\Omega)(X_1,\ldots,X_{m})U_{\alpha}
 -\sm{\sum}_{k}\cos\theta\,  B(X_k,X_k) \\
&& +\sm{\sum}_k\sm{\sum}_{i<k}(-1)^{k+i}(\Omega^{\sharp}(
B(X_k,X_i),X_1,\ldots, \hat{X_i},\ldots
,\hat{X_k},\ldots, X_m))^{\bot}\\
&&+\sm{\sum}_k 
\sm{\sum}_{k<i}(-1)^{k+i-1}(\Omega^{\sharp}(
B(X_k,X_i),X_1,\ldots, \hat{X_k},\ldots
,\hat{X_i},\ldots, X_m))^{\bot}.
\end{eqnarray*}
Interchanging $i$ by $k$ in the later line and using the
symmetry of $B$ and that $d\Omega=0$ we prove the first equality of the lemma. 
The computation of the second equality is similar with 
$\langle \nabla\Phi, B\rangle =$ 
$
\sum_{jk}\bar{g}(\nabla_{X_j}\Phi(X_k),B(X_j,X_k))$, 
recalling the definition in (1.6), where
$\langle \Psi, \wedge^2 B\rangle = \sum_{i<k}\sum_j\langle
\Psi(X_i,X_k), B(X_j,X_i)\wedge B(X_j,X_k)\rangle$.
\qed\\[5mm]
\noindent
{\bf \em  Proof of Theorem 1.2. ~}
Consider the vector field $Z$ on $M$ defined by
\begin{equation}
 g(Z,X)=\bar{g}(\Phi(X),H)~~~\forall X\in TM.
\end{equation}
 Using lemma 2.2
we have
\begin{eqnarray}
div(Z)&=&-\bar{g}( \delta \Phi, H)+ \sm{\sum}_i\bar{g}(\Phi(X_i),\nabla^{\bot}
_{X_i}H)\nonumber \\
&=& -m\cos\theta \, \|H\|^2 + \langle \nabla^{\bot}H, \Phi\rangle
+\bar{\nabla}_{H}\Omega(X_1,\ldots, X_m).
\end{eqnarray}
Assume $Z\neq 0$. Take $X={Z}/{\|Z\|}$ and
  $U={\Phi(X)}/{\|\Phi(X)\|}$. By  lemma 2.1,
$
\|\Phi(X)\|^2$ $
=\bar{g}(\Phi(X),U)^2\leq \sin^2\theta$.
From this inequality and applying  Schwartz inequality in (2.1) 
we get
\begin{equation}
\|Z\|\leq \sin\theta\, \|H\|.
\end{equation}
If $\nu$ denotes the outward unit of  $\partial D$,
integration of (2.2)  on  $D$ 
 and   (2.3) gives
\begin{eqnarray*}
\lefteqn{\left|\,  \int_D(-m\cos\theta \|H\|^2 
+\langle \nabla^{\bot}H,\Phi\rangle)dV
+\int_D \bar{\nabla}_{H}\Omega\, \right|} \\
&\leq& \left|\int_{\partial D}\langle Z,\nu\rangle dA\right|
\leq \int_{\partial D} \sin\theta\, \|H\|dA~~~~~~~\qed
\end{eqnarray*}
Theorem 1.3 and its corollary 1.1  are an  immediate consequence
 of theorem 1.2.\\[2mm]
{\bf \em Proof of Corollary 1.2. ~}
By theorem 1.2, 
$ m\cos\theta\int_{M}\|H\|^2dV=\int_M\langle
\nabla^{\bot}H,\Phi\rangle dV.$
We use Schwartz 
and a geometric-arithmetic inequality to obtain 
$$|\langle \nabla^{\bot}H,\Phi\rangle |\leq 
\sm{\sum}_i\left|\bar{g}(\nabla^{\bot}_{X_i}H,\Phi(X_i))\right|\leq 
\sm{\sum}_i\|\nabla^{\bot}_{X_i}H\|\, \mu\, \sin\theta\leq \sqrt{m}\, 
\mu\, \sin\theta\,\|\nabla^{\bot}H\|,$$
where $X_i$ is any orthonormal basis of $T_pM$.
If equality holds, then $\Phi(X_i)=\alpha_i\nabla^{\bot}_{X_i}H$
or $\nabla^{\bot}_{X_i}H=\beta_i\Phi(X_i)$, and if 
$\nabla^{\bot}_{X_i}H\neq 0$ we must have 
$\|\Phi(X_i)\|=\mu \,\sin\theta$,
where $\alpha_i,\beta_i\in 
\mathbb{R}$. Since $\mu\sin\theta\neq 0$ we must 
have $\nabla^{\bot}_{X_i}H=\Phi(\psi(X_i))$, $\forall i$.
\qed
\section{$\Delta \cos\theta$}
\renewcommand{\thesection}{\arabic{section}}
\renewcommand{\theequation}{\thesection.\arabic{equation}}
\setcounter{equation}{0}
In this section we are assuming $\Omega$ is parallel and 
$F:M^m\to \oM^{m+n}$ is an immersion with $\cos\theta>0$.
We use the curvature sign convention $\bar{R}(X,Y)=-[\bar{\nabla}_X,
\bar{\nabla}_Y]+\bar{\nabla}_{[X,Y]}$. Thus,
$\bar{R}(X,Y,Z,W)=\bar{g}(\bar{R}(X,Y)Z,W)$.
\begin{lemma} $ \nabla\cos\theta= B_{\Phi}^{\sharp}$ ~and 
$$\begin{array}{rcl}
 \Delta\cos\theta &=&-\cos\theta\, \tilde{Q}(B)
+{m}\langle \nabla^{\bot}H,\Phi\rangle
-\sm{\sum}_{ij}\bar{R}(X_i,X_j,X_i,\Phi(X_j))\\
\Delta\log(\cos\theta) &=&-Q_{\Omega}(B) + \frac{m}
{\cos\theta}\langle \nabla^{\bot}H,\Phi\rangle
-\frac{1}{\cos\theta}\sm{\sum}_{ij}\bar{R}(X_i,X_j,X_i,\Phi(X_j))\\
\Delta\La{(}\frac{1}{\cos\theta}\La{)} &=& \frac{1}{\cos\theta}\La{(}
\hat{Q}_{\Omega}(B)- \frac{m}
{\cos\theta}\langle \nabla^{\bot}H,\Phi\rangle
+\frac{1}{\cos\theta}\sm{\sum}_{ij}\bar{R}(X_i,X_j,X_i,\Phi(X_j))\La{)}
\end{array}$$
where $\hat{Q}_{\Omega}(B)
=Q_{\Omega}(B)+\frac{1}{\cos^2\theta}\|B_{\Phi}\|^2$.
Moreover, $\|\nabla\cos\theta\|^2\leq m\, \sin^2\theta\|B\|^2$
and $\|\Phi\|^2\leq m\, \sin^2\theta$. 
Furthermore, if $n=1$, $\|\Phi\|^2\geq \sin^2\theta
$. In the later case
if $m=2$  and $H=0$, then 
$\|\nabla \cos\theta\|^2=\frac{1}{2}
\|B\|^2(\|\Omega\|^2-\cos^2\theta)=\frac{1}{2}
\|B\|^2\|\Phi\|^2$.
\end{lemma}
\noindent
\em Remark. \em In this lemma 
the expression of $\Delta\cos\theta$ is still valid for $\cos\theta$ with
any value in $\mathbb{R}$, since, at points where $\cos\theta=0$,
$\cos\theta\tilde{Q}_{\Omega}(B)$ means
$2\langle \Psi,\wedge^2B\rangle$.\\[2mm]
\em Proof. \em  For a local orthonormal frame $X_i$,
and where $(j)$ denotes "place $j$", 
\begin{eqnarray*}
\lefteqn{d\cos\theta (X_k) = \sm{\sum_j} \Omega(X_1,\ldots,
\bar{\nabla}_{X_k}{X_j}_{\ti{(j)}},\ldots, X_m)} \\
&=& \sm{\sum_j} \Omega(X_1,\ldots, B(X_k, X_j)_{\ti{(j)}},
\ldots, X_m)
+\sm{\sum}_j \Omega(X_1,\ldots,  \lnab{X_k}{ X_j}_{\ti{(j)}},
\ldots, X_m)\\
&=& \sm{\sum_j} \Omega(X_1,\ldots, B(X_k, X_j)_{\ti{(j)}},
\ldots, X_m)=B_{\Phi}(X_k), 
\end{eqnarray*}
where we used in the last equality
$\langle \lnab{X_k}X_j, X_j\rangle =0$,
 and
$\Omega(X_1,\ldots,  { X_k}_{\ti{(j)}}
\ldots, X_m)$ vanish  for $k\neq j$.
Differentiation  of the later equation at $p$, and using that 
\begin{eqnarray*}
\bar{\nabla}_{X_k}B(X_k,X_j)&=&{\nabla}_{X_k}B(X_k,X_j)-
\sm{\sum}_i \bar{g} (B(X_k,X_j),B(X_k,X_i))X_i\\
\sm{\sum_{k}}\nabla_{X_k}B(X_k,X_j)
&=&\sm{\sum_{k}}\nabla_{X_k}B(X_j,X_k)
= m\lnabo{X_j}H -\sm{ \sum_{k}}(\bar{R}(X_k,X_j) X_k)^{\bot},
\end{eqnarray*}
where in the latter equality we used the
  Codazzi's equation for $B$, 
$\nabla_{X_k}B(X_j,X_k)= $ $\nabla_{X_j}B(X_k,X_k)
- (\bar{R}(X_k,X_j) X_k)^{\bot}$,
we have
\begin{eqnarray*}
\triangle \cos\theta &=&
\sm{\sum_{k}}\lnab{X_{k}}d\cos\theta (X_{k})\nonumber \\
&=& \sm{\sum_k \sum_{s<j}}2\Omega(X_1, \ldots,
 B(X_s,X_k)_{\ti{(s)}},\ldots,
B(X_j,X_k)_{\ti{(j)}}, \ldots, X_m)\nonumber \\
&&+\sm{\sum_{kj}} \Omega(X_1,\ldots,\nabla_{X_k}B(X_k,X_j)
_{\ti{(j)}}, \ldots, X_m)
-\cos\theta \bar{g}( B(X_k,X_j),
 B(X_k, X_j)) \nonumber\\
&=& 
\sm{\sum}_{s<j}\sm{\sum}_k 2\langle \Psi(X_s,X_j),B(X_s,X_k) \wedge
B(X_j,X_k)\rangle + m\langle \nabla^{\bot}H,\Phi\rangle\\
&& -\cos\theta \|B\|^2
-\sm{\sum}_{kj} \bar{R}(X_k,X_j,X_k,\Phi(X_j)).
\end{eqnarray*}
The lemma now follows from the expressions of
$\Delta \log \cos\theta$ and $\Delta (\frac{1}{\cos\theta})$
in terms of $\Delta\cos\theta$.
Next we estimate $\|\nabla \cos\theta\|^2$ and $\|\Phi\|$.
By lemma 2.1, $\|\Phi(X_k)\| \leq \sin\theta$, and so
$\|\Phi\|^2 \leq m\, \sin^2\theta$. Now, from the first equation 
in this proof,
\begin{eqnarray*}
\|\nabla\cos\theta\|^2 &=&\|B_{\Phi}\|^2
\leq \sm{\sum}_{ijk}\sin^2\theta\|B_{ij}\|\, \|B_{ik}\| \nonumber\\
&\leq & \sm{\sum} _{ijk}\frac{\sin^2\theta}{2}
(\|B_{ij}\|^2 + \|B_{ik}\|^2)\leq m\sin^2\theta\, \|B\|^2,
\end{eqnarray*}
where $B_{ij}=B(X_i,X_j)$.
Note that
$\|\Omega\|^2\geq 1$, because there exists
a calibrated subspace. Consequently, if  $n=1$, 
$1\leq \|\Omega\|^2=\cos^2\theta +\|\Phi\|^2$ what implies
$\|\Phi\|^2\geq \sin^2\theta$.
Finally in case $m=2$ and $n=1$, let $\nu$ be a unit normal to $M$ and
set $B_{ij}=\bar{g}(B(X_i,X_j),\nu)$.
From minimality of $F$
\begin{eqnarray*}
\|\nabla\cos\theta\|^2 
&=& \Omega(B(X_1,X_1), X_2)^2+
\Omega(B(X_2,X_1), X_1)^2\\
&&+
\Omega(B(X_1,X_2), X_2)^2+
\Omega(B(X_2,X_2), X_1)^2\\
&=&B_{11}^2\Omega(\nu,X_2)^2+B_{12}^2\Omega(\nu,X_1)^2
+B_{12}^2\Omega(\nu,X_2)^2 +B_{22}^2\Omega(\nu,X_1)^2\\
&=& \frac{1}{2}\|B\|^2(\Omega(\nu,X_2)^2+\Omega(\nu,X_1)^2)
=\frac{1}{2}\|B\|^2(\|\Omega\|^2-\Omega(X_1,X_2)^2).\qed
\end{eqnarray*}
An immediate consequence  from 
the last equality of lemma 3.1 follows
next:
\begin{proposition}  If $m=2$ and $n=1$, $F$ is minimal
 and $\cos\theta $ is constant, then
either $F$ is a calibrated submanifold or it is totally geodesic.
\end{proposition}
\begin{proposition} If $F$ is a parallel submanifold then $(1.7)=0$ holds.
\end{proposition}
\noindent
{\em  Proof. \em } From proof of lemma 3.1
$\sum_k\nabla_{X_k}B(X_k,X_j)=m\nabla_{X_j}H
-\sum_k(\bar{R}(X_k,X_j)X_k)^{\bot}.$
Since $\nabla B=0$, then  $\sum_k(\bar{R}(X_k,X_j)X_k)^{\bot}=0.$\qed
\\[4mm]
{\bf \em Proof  of Proposition 1.3 } 
Since $F$ is calibrated, $F$ is minimal,
and  by lemma 2.1 $\Phi=0$. From lemma 3.1,  
$0=Q_{\Omega}(B)=\tilde{Q}_{\Omega}(B)$, and so
$B=0$.\qed\\[4mm]
Recall the average value of a function $f$ on
a domain $D$ and on $\partial D$ is given by:
$$-\!\!\!\!\!\!\!\!\!\int_D f dV= \frac{1}{V(D)}\int_D f dV, ~~~~~~~~
-\!\!\!\!\!\!\!\!\!\int_{\partial D} f dA= \frac{1}{A(\partial D)}
\int_{\partial D} fdA.\\[1mm]$$
{\bf \em Proof of Theorem 1.5.~ } 
Using lemma 3.1 we obtain the first inequality of (1.8), and
under the assumptions of the 
theorem  have
\begin{equation}
 \Delta \log \cos\theta \leq -Q_{\Omega}(B)\leq -\delta\|B\|^2.
\end{equation}
Thus, applying lemmas 2.1 and 3.1,  
 we have after integration of (3.1),
$$\delta \int_D\|B\|^2\leq \int_{\partial D} -g(\nabla \log\cos\theta,
\nu)dV\leq \int_{\partial D}\sqrt{m}\sup_{\partial D}\tan\theta \|B\|dA$$
and second inequality of  (1.8) is proved. If $\tan\theta \|B\|$ is integrable
on $M$, then so it is $\nabla \log\cos\theta$, by the first inequality
of (1.8). Since
$-\log\cos\theta$ is an subharmonic function by (3.1),
      then applying the Stokes theorem for complete manifolds in the version
  given by Yau (\cite{yau} Corollary page 660), we conclude
$\Delta\log\cos\theta=0$.  From (3.1), this implies 
$B=0$  on $
D$, and (A) is proved.  
(B) follows immediately.
To prove (C) we use the Omori-Cheng-Yau maximum principle.
Under the assumptions, 
 by Gauss equation $Ricci^M$ is bounded from below. If we
 assume $\inf_M\cos\theta>0$,
then we take a sequence $p_k$ such that $u(p_k)\to \sup_M u$,
$\nabla u(p_k)\to 0$ and $\lim_k\Delta u(p_k)\leq 0$ when
$k\to +\infty$, where $u= \log(\cos\theta)^{-1}$. This implies by lemma
3.1 that ${Q}_{\Omega}(B)(x_k)\to 0$, and so $\inf_M\|B\|=0$.
\qed
\\[5mm]
{\bf \em Proof of Theorem 1.6.~ } 
Recall that if a surface is parabolic,  any nonnegative superharmonic
function is constant.
By Gauss equation, 
the sectional curvature of $M$ satisfies
$ K_M= \bar{R}(X_1,X_2,X_1,X_2)-\frac{1}{2}\|B\|^2.$
Let ${C}$ be a compact set of $M$ such that
$\bar{K}\circ F\geq 0$ away from ${C}$.
(1) Since $\oM$ is 3-dimensional,
  $Q_{\Omega}(B)\geq \|B\|^2$  and 
by lemma 3.1,
$\|\nabla \cos\theta\|^2\geq \frac{1}{2}\|B\|^2\sin^2\theta$. 
Again, by lemma 3.1
and the assumptions of the theorem we have
$\Delta\cos\theta \leq -\cos\theta\|B\|^2$ and so
\begin{eqnarray}
\Delta\log(1+\cos\theta) &=& \frac{1}{(1+\cos\theta)}\La{(}\Delta
\cos\theta-\frac{\|\nabla\cos\theta\|^2}{(1+\cos\theta)}\La{)} \nonumber
\\
&\leq& \frac{1}{(1+\cos\theta)}\La{(}-\cos\theta \|B\|^2
-\frac{\|\nabla\cos\theta\|^2}{(1+\cos\theta)}\La{)}\leq 0
\end{eqnarray}
We consider on $M$ the complete
metric $\tilde{g}=(1+\cos\theta)^pg$, where we choose
$p\geq 2$.
This metric
has sectional curvature $\tilde{K}$ that satisfies on $M\sim C$
\begin{eqnarray*}
\tilde{K} &=& \frac{1}{(1+\cos\theta)^p}\la{(}\bar{R}(X_1,X_2,X_1,X_2)
-\frac{1}{2}\|B\|^2-\frac{p}{2}\Delta\log(1+\cos\theta)\la{)}\\
&\geq & \frac{1}{(1+\cos\theta)^p}\La{(}
-\frac{1}{2}\|B\|^2+ \frac{p }{2}\frac{1}{(1+\cos\theta)}
(\cos\theta \|B\|^2 +\frac{\|\nabla \cos\theta\|^2}
{(1+\cos\theta)})\La{)}\\
&\geq&\frac{1}{(1+\cos\theta)^p}\La{(}
 -\frac{1}{2} +\frac{ p }{2}\frac{\cos\theta}{(1+\cos\theta)}
+\frac{p }{4}\frac{\sin^2\theta}{(1+\cos\theta)^2}\La{)}\|B\|^2\\
&= & \frac{p-2}{4(1+\cos\theta)^{p}}\|B\|^2.
\end{eqnarray*}
Note that $d\tilde{V}=(1+\cos\theta)^p dV$ and so
$\int_M \tilde{K}^-d\tilde{V}<+\infty$, where $\tilde{K}^-=
\max\{-\tilde{K},0\}$. This implies that $(M,\tilde{g})$ is
parabolic, and  so it is $(M,g)$, since
 $\tilde{\Delta}=(1+\cos\theta)^{-p} \Delta$.
From (3.2) 
 we have $\cos\theta$ constant and
since $\cos\theta>0$,   we conclude  
 that $B=0$.
(2)
We consider on $M$ the metric $\tilde{g}=\cos^{\frac{1}{\delta}}\theta g$.
This metric is complete because
$1\geq \cos\theta\geq \tau>0$,
and by lemma 3.1 and the assumptions on the theorem,
 the sectional curvature  satisfies on
$M\sim C$,
$$\tilde{K} = {\cos^{-\frac{1}{\delta}}\theta}
(K_M-\frac{1}{2\delta}\Delta \log \cos\theta)
\geq {\cos^{-\frac{1}{\delta}}\theta} \bar{R}(X_1,X_2,X_1,X_2)
\geq 0$$
and 
$\Delta\left(\frac{1}{\cos\theta}-\frac{1}{\tau}\right)\geq \frac{\delta}{
\cos\theta}\|B\|^2$.
By the same arguments as in (1), $B=0$.
\qed\\[5mm] 
{\bf \em  Proof of corollaries 1.3 and 1.4. } 
First we note that the graph $\Gamma_f$ of a map $f:(M,g_1)\to (N,h)$
defines a complete submanifold of
$(M\times N, g_1\times h)$ provided $(M,g_1)$ is
complete.  
We are considering $\Omega$ the volume calibration.
(1) is a result of theorem 1.6(1), since  graphs
satisfy $\cos\theta>0$. (2) is a consequence  of theorem 1.6(2), 
because  the eigenvalues of $f^*h$ satisfy
$|\lambda_1\lambda_2|\leq 1-\delta$, what  implies $Q_{\Omega}(B)\geq 
\delta\|B\|^2$ and  $\cos\theta
\geq \tau= (1+C+(1-\delta)^2)^{-1/2}$, where  $\|df\|^2\leq C$
(see subsection 5.1).
 To prove (3) we use 
theorem 1.5 (A). Since $\lambda_i^2\leq 1-\delta$ (case $n\geq 2$),
$Q_{\Omega}$ is
$\delta$-positive. Now
we only have to check that (1.7) holds.
Using $X_i$ and $X_{m+\alpha}$ suitable orthonormal
frames of $(TM,g=g_1+f^*h)$ and of $(NM,\bar{g})$ respectively 
(see  section 5),
we have
\begin{eqnarray}
\lefteqn{\sm{\sum_{i\neq j}\bar{R}(X_i,X_j,X_i,\Phi(X_j))=}} \nonumber
\\
&=&\!\!\!\sm{\cos\theta
\sum_{i\neq j}\sm{\frac{\lambda_j^2}{(1+\lambda_i^2)(1+\lambda_j^2)}}\left(
R_{1}^M(a_i,a_j,a_i,a_j)-\lambda_{i}^2
R^N(a_{i+m},a_{j+m},a_{i+m},a_{j+m})\right)}~~~~~~\\
&=&\!\!\!\sm{ \cos\theta\sum_{j}\sum_{i\neq j}\La{(}\sm{\frac{\lambda_j^2}
{(1+\lambda_j^2)}}
 K_1(a_i,a_j) - \sm{\frac{\lambda_i^2\lambda_j^2}
{(1+\lambda_i^2)(1+\lambda_j^2)}}(K_1(a_i,a_j)+K_N(a_{i+m},a_{j+m}))\La{)}}
\nonumber\\
&=&\!\!\!\sm{\cos\theta\La{(}\sum_{j}\sm{\frac{\lambda_j^2}{(1+\lambda_j^2)}}
 Ricci_{1}(a_ja_j) - \sum_{j\neq i}\sm{\frac{\lambda_i^2\lambda_j^2}
{(1+\lambda_i^2)(1+\lambda_j^2)}}(K_1(a_i,a_j)+K_N(a_{i+m},a_{j+m}))\!
\La{)}}~~~~~~
\end{eqnarray}
If $n=1$, $\lambda_i=0$ for $i\geq 2$ and the last term of (3.4) disappears.
For $n\geq 2$, if we assume at each point $p\in M$ and two-planes 
$P$ of $T_pM$
and   $P'$ of $T_{f(p)}N$
$K_1(P)\geq K_N(P')^+=\max\{K_N(P'), 0\}$ (
or $Ricci_{1}(p)\geq 0$ and, $K_N(P')\leq -K_1(P)$ respectively )
then  we get 
$\sum_{i\neq j}\bar{R}(X_i,X_j,X_i,\Phi(X_j))\geq 0$ from (3.3)
(from (3.4) respectively). This implies by theorem 1.5,
that $\Gamma_f$ is a totally geodesic
submanifold of $\oM$, on therefore
 $f:(M,g_1)\to (N,h)$ is also totally geodesic
( see \cite{[S2]} or \cite{[lisal],[S3]}),
and so we have equality to zero in all above equalities (see proposition
3.2). In this case all $\lambda_i$ are constant.
Moreover if at some point
$K_1(p)>0$ ( or $Ricci_{1}>0$ respectively) then $\lambda_j=0$ for all $j$,
and $f$ is constant. \qed\\[3mm]
{\bf Proof of corollary 1.5.} We have
$\cos\theta>0$, $\tilde{Q}_{\Omega}(B)={Q}_{\Omega}(B)\geq\|B\|^2$ and
(1.7) holds as well. In this case we use a modified version of the proof of
theorem 1.5 (A), by considering $\Delta \cos\theta=-\cos\theta
\tilde{Q}_{\Omega}(B)\leq 0$. 
Integrability of $\|\nabla \cos\theta\|\leq \sin\theta\|B\|$ implies  
$\Delta \cos\theta=0$, and
so $\|B\|=0$, that is $M$ is an hyperplane.\qed

\section{The Cheeger constant of a submanifold}
\renewcommand{\thesection}{\arabic{section}}
\renewcommand{\theequation}{\thesection.\arabic{equation}}
\setcounter{equation}{0}
In this section we estimate the Cheeger constant (1.1) of
a Riemannian manifold
$(M,g)$. If $\partial M\neq \emptyset$ then $D$ satisfies $\partial D\cap
\partial M=\emptyset$. If $M$ closed 
we may let $D=M$.
\begin{proposition} \cite{yau2} If $M$ is complete
simply connected and the sectional curvature  satisfies
$K_M\leq -K$, $K$ a positive constant, then $\h(M)\geq (m-1)\sqrt{K}$.
\end{proposition}
\noindent
The proof is based on the use of the comparison theorem
to obtain
 $\Delta r\geq (m-1)\sqrt{K}$ where $r$ is the distance function
to a point, and 
integration on a domain $D$.   

\begin{proposition} If $M$ is complete and $Ricci^M\geq 0$ then 
$\h(M)=0$.
\end{proposition}
\noindent
\em Proof. \em 
If  we assume the Ricci curvature of $M$ satisfies
$Ricci^M\geq 0$, following \cite{[ADR]} for $m=2$,
by a result due to Cheng \cite{[Cheng]}
the first eigenvalue of the Dirichlet problem
on a geodesic ball $B_r(p)$ is less than or equal to the first
eigenvalue of a geodesic ball of the same radius of $\R^m$,
that is $C_1/r^2$ for some constant $C_1>0$ that does not depend on $r$.
Therefore
$\lambda_1(B_r(p))\leq {C_1}/{r^2}$, for $ 0<r<+ \infty$. 
By a well known inequality
due to Cheeger (Theorem 3 p.95 in \cite{[Cha]}), we get
$
\h^2(M)\leq \h^2(B_r(p))\leq 4 \lambda_1(B_r(p))\leq \frac{4C_1}{r^2}.
$
 This  implies for $M$ complete that
$ \h (M)=0$.\qed\\[5mm]
Given a   smooth function
$f:M\to \mathbb{R}$ and a regular value of $|f|$,
$t\in R_{|f|}$, 
 the  sets 
$$\begin{array}{c} D^+_{f}(t)=\{ p\in M: |f(p)|\geq t\},~~~~~
D^-_{f}(t)=\{ p\in M: |f(p)|\leq t\},\\[2mm]
\Sigma_{f}(t)=\{ p\in M: |f(p)|=t\}
\end{array}$$
define smooth submanifolds 
with $\partial D^{\pm}_f(t)=\Sigma_f(t)$.
Set $V_{\pm}(t)=V(D^{\pm}_{f}(t))$ and $A(t)=A(\Sigma_{f}(t))$.
Then $V_{\pm}(t)$ are smooth on $R_{|f|}$ and 
the co-area formula (see e.g. \cite{[Cha]}) 
states that for any nonnegative mensurable function $h$ (or $h\in L^1(M)$),
\begin{eqnarray*}
\int_Mh\|\nabla f\|dV &=&\int_{0}
^{+\infty}dt\int_{\Sigma_f(t)}h\,  dA(t). 
\end{eqnarray*}
Applying the co-area formula to $h=\|\nabla f\|^{-1}\mathcal{X}_{A}$,
with $A=D_f^{\mp}(s)$
for $s$ regular value, where $\mathcal{X}_A$ is the 
characteristic function of a set $A$, one obtains
\begin{eqnarray*}
{V_{\mp}}'(s) = \pm\int_{\Sigma_{f}(s)}\|\nabla f\|^{-1}dA(s).
\end{eqnarray*}
\begin{lemma} If $D^{\mp}_f(s)$ 
 is  compact for  $s< \sup |f|$ ~($
s> \inf |f|$ resp.), then
$$\h(M)\leq 
\pm \frac{-\!\!\!\!\!\!\!\int_{D^{\mp}_f(s)}\|\nabla f\|dV}{ 
s--\!\!\!\!\!\!\!\int_{D^{\mp}_f(s)}|f|dV}.
$$
\end{lemma}
\noindent
\em Proof. \em
Using the co-area formulas for $f$ restricted to the interior of
${D_f^-(s)}$ ( at regular values $s$), 
 and that $\h(M)\leq A(t)/V_{-}(t)$, 
$\forall t<s$,  we have
\begin{eqnarray*}
\int_{D_f^-(s)}\|\nabla f\|dV&=& \int_0^{s} A(t)dt\geq \int_0^{s} 
\h(M)V_-(t)dt\\
&=&\h(M)\LA{(} tV_-(t)\LA{]}^{s}_0 -\int_0^{s}tV'_-(t)\LA{)}\\
& =& \h(M)\LA{(}sV_-(s)-\int_0^{s}\int_{\Sigma_f(t)}
|f(x)|\|\nabla f\|^{-1}dA(t)dt\LA{)}\\
&=& 
\h(M)\La{(}sV_-(s)-\int_{D_f^-(s)}|f(x)|dV\LA{)}
\end{eqnarray*}
where in the  last equality we use the co-area formula for
$h=\frac{|f|}{\|\nabla f\|}$. 
Similarly for 
$D_f^+(s)$. Note that  the functions
 $(s--\!\!\!\!\!\!\!\int_{D^{-}_f(s)}|f|dV)=\int_0^sV_{-}(t)dt$
and $(-s+-\!\!\!\!\!\!\!\int_{D^{+}_f(s)}|f|dV)=\int_s^{t_+}V_{+}(t)dt$,
where $t_+=\sup_M|f|$,
are increasing  on $s$.\qed
\begin{corollary} The Cheeger constant of $M$ vanish if 
there exist a smooth function $f$ such that
 $D^-_f(s)$ is compact  $\forall s<\sup_M|f|=+\infty$, 
and for some constants $0<\alpha, \delta <1$, we have
$-\!\!\!\!\!\!\!\int_{D^{-}_f(s)}\|\nabla f\|dV\leq s^{\alpha}$, and 
$-\!\!\!\!\!\!\!\int_{D_f^-(s)}|f|\leq (1-\delta)s,$ when $s\to +\infty$.
\end{corollary}
\noindent
Let $r$ be the distance function to a point $p$. 
Recalling that $\|\nabla r\|=1$, we have
\begin{corollary}  If $r^2$ is smooth for $r<s$, ~ 
$\h(M)\leq -\!\!\!\!\!\!\!\int_{B_{s}(p)}2krdV
/\la{(}
{s^2- -\!\!\!\!\!\!\!\int_{B_{s}(p)}r^{2}dV}\la{)}$.
In particular, if $-\!\!\!\!\!\!\!\int_{B_{s}(p)}
r^{2}dV\leq (1-\delta) s^{2}$, 
where $0< \delta<1$ is a constant, 
 $(M,g)$ is complete, and $r^2$ smooth on $M$,  then $\h(M)=0$.
\end{corollary}
\noindent
Note that in the previous corollary, $-\!\!\!\!\!\!\!\int_{B_{s}(p)}
rdV\leq (-\!\!\!\!\!\!\!\int_{B_{s}(p)}
r^{2}dV)^{1/2}\leq \sqrt{(1-\delta)} s$, and $D_f^-(s^2)=B_{s}(p)$.
Thus, if $\h(M)\neq 0$, then $\lim_{s\to +\infty}
\frac{1}{s^2}{-\!\!\!\!\!\!\!\int_{B_{s}(p)}}r^{2}dV=1$.
If $M=\mathbb{R}^m$, in corollary 4.2 we may take $\delta=\frac{m}{2+m}$.
\\[3mm]
If $M$ is a  submanifold of a $(m+n)$-dimensional
Riemannian manifold $(\oM,\bar{g})$ we can also estimate 
the Cheeger constant of $M$ under certain conditions.
Recall (\cite{[S1],[RS]}) that a vector field $\bar{X}$ on $\oM$
is  \em strongly convex \em on a open set $U$ of $\oM$ if 
$$L_{\bar{X}}\bar{g}\geq 2\alpha \bar{g}$$
where $\alpha>0$ is a constant. Examples of such vector fields
are $\sm{\frac{1}{2}}\bar{\nabla} \bar{r}^2=r\frac{\partial}{
\partial r}$ on a
geodesic ball of $\oM$ of radius $R$ and center $\bar{p}$
   that does not intercept the cut locus at $\bar{p}$ and
$\sqrt{\kappa}R<\pi/2$ where $\kappa=\max\{0, \sup_{B_R(\bar p)}\bar{K}\}$
and $\bar{K}$ are the sectional curvatures of $\oM$. A strictly 
convex function $f$
on $\oM$ with $Hess \, f\geq \alpha \bar{g}$ defines a strongly convex
vector field $\nabla f$.
Positive homothetic no-Killing vector fields  are  strongly convex.
In $\mathbb{R}^{m+n}$ the position vector field
$\bar{X}_x=x$ is such an example.
If $F:M\to \oM$ is an immersed submanifold let
$X_F$ denote the vector field $X$ along $F$.

\begin{lemma} If $\oM$ carries a strongly convex vector field
$\bar X$ on a neighbourhood of an immersion
$F:M\to \oM$ then  
$(\sup_M\|{\bar X}_F\|)^{-1}\leq \frac{1}{\alpha}(
\frac{1}{m}\h(M)+\sup_M\|H\|)$.
In particular if $M$ is minimal and 
 has zero Cheeger constant then ${\bar X}_F$ is
unbounded.
\end{lemma}
\noindent
\em Proof. \em  By an elementary computation
(see \cite{[HoffSpru]} or \cite{[S1]})
for any immersion $F$
$$m\,\bar{g}(H,{\bar X}_F) = div_g({\bar X}^T_F) 
-\sm{\frac{1}{2}}\, tr_g \, L_{\bar X}\bar{g}$$
where ${\bar X}_F^T$ is the projection of $\bar{X}_F$ onto
$TM$. Integration on a domain $D$
gives
\begin{eqnarray*}
\alpha mV(D) & \leq & \int_{\partial D}g({\bar X}^T_F,\nu)dA
-\int_Dm\,\bar{g}(H,{\bar X}_F)dV\\
&\leq& \sup_D\|{\bar X}_F\|\La{(}A(\partial D)+ \int_Dm\,\|H\|dV\La{)}
\end{eqnarray*}
where  $\nu$ is a unit normal to $\partial D$. Thus
$(\sup_D\|{\bar X}_F\|)^{-1}\leq \frac{1}{\alpha }\La{(}
\frac{1}{m}
\frac{A(\partial D)}{V(D)}+-\!\!\!\!\!\!\!\!\int_D\|H\|dV \La{)}$,
with $-\!\!\!\!\!\!\!\!\int_D\|H\|dV\leq sup_M\|H\|$.
Taking the infimum on $D$ we obtain the proposition.\qed
\begin{corollary} If $F:M\to \mathbb{R}^{n+m}$ is a minimal immersion
with zero Cheeger constant, then $F$ is unbounded.
\end{corollary}
\noindent
{ \bf \em  Proof of Theorem 1.4. }
Using theorem 1.3, on each ball $B_r(p)$, and for
any domain $D\subset B_r(p)$,
$\|H\|\leq \frac{1}{m}(\inf_{B_r(p)}\cos\theta)^{-1}
\frac{A(\partial D)}{V(D)}$,
and so, taking the infimum for $D\subset B_r(p)$,
$m\|H\|\leq (\inf_{B_r(p)}\cos\theta)^{-1}
\h(B_r(p))$.
By assumption of the theorem
$\inf_{B_r(p)}\cos\theta\geq Cr^{-\beta}$, and (1.4)  leads to
$\|H\|\leq C {r^{\beta-1}}$
for some constant $C>0$ that does not depend on $r$. Thus, letting
$r\ra +\infty$ we obtain $\|H\|=0$, and theorem 1.4 is proved.\qed \\[4mm]
\noindent
Proposition 1.1 is an immediate
consequence of  proposition 4.2 and Theorem 1.1.\\

If under certain conditions we have $\Delta\cos\theta\leq 0$,
(as  in theorem 1.5) and $\theta$ is not constant,  by the maximum principle,
for any regular value of $\cos\theta$,  $D^-(\epsilon)=\{p:\cos\theta
\leq \epsilon\} $ cannot be a compact domain.  Next we assume for
 $\epsilon > 0$, the set
$D^+({\epsilon})=\{p: \cos\theta \geq \epsilon\} $ to be compact.
The next proposition is  an attempt to understand what happens if one 
replaces the assumption $Q(B) \geq  \delta \|B\|^2$  by the weaker condition 
$Q(B)\geq  0.$ .
\begin{proposition} Assume $\cos\theta>0$, and $D^+({\epsilon})$ is compact
$\forall \epsilon \in (0,1]$. Then:\\[2mm]
(a)~If there exists  constants $\alpha, \delta\in (0,1)$ such that
$-\!\!\!\!\!\!\!\int_{D^+({\epsilon})}(\cos\theta)^{-1}
dV$ $\leq (1-\delta)\frac{1}{\epsilon}$, 
and  $ -\!\!\!\!\!\!\!\int_{D^+({\epsilon})}
\sin\theta\cos\theta^{-2}\|B\|\leq (\frac{1}{\epsilon})^{\alpha}$, 
 for
 $\epsilon\to 0$,
 then $\h(M)=0$.\\[2mm]
(b) If $M$ is immersed with parallel mean curvature,
 ${Q}_{\Omega}(B) +(\cos\theta)^{-2}\|B_{\Phi}\|^2
\geq  0$  and (1.7) holds,  and for some  constant $\alpha>1$,
$\int_M {(\cos\theta)^{-(\alpha+2)}}
\sin\theta\,\|\bar{R}\|dV
<+\infty$,
 then either
 $\int_M{(\cos\theta)^{-(\alpha+4)}} \|B\|^2dV=+\infty,$
or $M$ is compact.
\end{proposition}
\noindent
{\bf \em Proof. }
By lemma 3.1, 
$f=1/\cos\theta$ in (a) satisfies  $-\!\!\!\!\!\!\!\int_{D^+({\epsilon})}
\|\nabla f\|\leq (\frac{1}{\epsilon})^{\alpha}$, and in (b)
$\Delta f\geq  0$ and $f>0$. Then (a) follows from  corollary 4.1.
Now we prove (b).
For each $s$ fixed we consider the compact sets
$ D_f^-(s)$ and $\Sigma_f(s)$  and  follow close the proof
in  \cite{PeterLi} (lemma 7.1)  replacing
$r$ by $f$. 
We take a cut off function 
$\phi(f):M\to \mathbb{R}_0^+$ where $\phi:\mathbb{R}_0^+\to [0,2]$ 
is a smooth nonnegative bounded
function  satisfying 
$\phi(t)=1$ if $t\leq s$, $\phi(t)=0$ if $t\geq 2s$, 
and $ (\phi')^2\leq Cs^{-2}$ and $|\phi''|\leq Cs^{-2}$,
$C>0$ a constant that does not depend on $s$
(see for example, lemmas 7.1 and 6.1 of \cite{PeterLi} and p. 661 of
\cite{yau}).
Assume $2s\in R_f$.
Then integrating
 $\Delta(\phi^2(f)f^{\alpha-1} f)$ on $D_f^-(2s)$, applying Stokes
and using that $\Delta f^{\alpha-1}\geq 0$ 
we have
\begin{eqnarray*}
\lefteqn{0\leq \int_{D_f^-(2s)}\phi^2(f)f^{\alpha -1}\Delta f \, dV}\\
&\leq& \int_{D_f^-(2s)}
 -\La{(}2\phi(f)\Delta(\phi(f))+2\|\nabla (\phi(f))\|^2\La{)} f^{\alpha}dV\\
&&\!-\int_{D_f^-(2s)}\!\!\left(4(\alpha-1)\phi(f)\langle \nabla (\phi(f)),
\nabla f\rangle  f^{\alpha-1} +2
\langle \nabla\La{(}\phi^2(f)f^{\alpha -1}\La{)},
\nabla f\rangle \right)dV.
\end{eqnarray*}
Hence,
\begin{eqnarray*}
\lefteqn{2\int_{D_f^-(2s)}\phi^2(f)(\alpha -1)f^{\alpha -2}\|\nabla f\|^2\leq}\\
&&\int_{D_f^-(2s)}-2\phi(f)\La{(}\phi''(f)\|\nabla f\|^2
+\phi'(f)\Delta f\La{)}f^{\alpha} dV\\
&&-\int_{D_f^-(2s)}2\left(\|\nabla(\phi(f))\|^2f^{\alpha}+
2(\alpha-1)\phi(f)\langle \nabla(\phi(f)),\nabla f\rangle f^{\alpha-1}
\right)dV\\
&& -4\int_{D_f^-(2s)}\phi(f)\langle \nabla(\phi(f)),
\nabla f\rangle f^{\alpha -1} \, dV\\
&\leq & \int_{D_f^-(2s)}\frac{\tilde{C}}{s^2}\|\nabla f\|^2f^{\alpha}
+\frac{\tilde{C}}{s}\Delta f \, f^{\alpha}
+(\alpha -1)\,\frac{\tilde{C}}{s}  f^{\alpha -1}\|\nabla f\|^2\\
&&+ \int_{D_f^-(2s)}-4\phi(f) f^{\alpha-1}
\langle \nabla (\phi(f)),\nabla f\rangle,
\end{eqnarray*}
where $\tilde{C}>0$ denotes a constant that does not depend on $s$.
On the other hand 
\begin{eqnarray*}
\lefteqn{-4\phi(f) f^{\alpha -1}\langle \nabla (\phi(f)), \nabla f\rangle
=-4\langle f^{\frac{\alpha}{2}}\nabla (\phi(f)),
f^{\frac{\alpha-2}{2}}\phi(f)\nabla f\rangle}\\
&\leq & 2\La{(}\frac{2}{(\alpha -1)}\|\nabla (\phi(f))\|^2f^{\alpha}
+\frac{(\alpha -1)}{2}f^{\alpha -2}\phi(f)^2\|\nabla f\|^2\La{)}.
\end{eqnarray*}
This implies
\begin{eqnarray*}
\lefteqn{\int_{D_f^-(s)}(\alpha -1)
f^{\alpha -2}\|\nabla f\|^2dV\leq
\int_{D_f^-(2s)}(\alpha -1)\phi^2(f)\,
f^{\alpha -2}\|\nabla f\|^2dV\leq }\\
&\leq& \int_{D_f^-(2s)}\La{(}\frac{\tilde{C}}{s^2}\|\nabla f\|^2f^{\alpha}
+\frac{\tilde{C}}{s}\Delta f \, f^{\alpha}\La{)}dV\\ 
&&+\int_{D_f^-(2s)}\La{(}(\alpha -1)f^{\alpha -1}\frac{\tilde{C}}{s}
\|\nabla f\|^2
+\frac{4}{(\alpha -1)}\frac{\tilde{C}}{s^2}\|\nabla f\|^2 f^{\alpha}
\La{)} dV.
\end{eqnarray*}
Assuming $\frac{\|B\|^2}{\cos^{4+\alpha}\theta}$ integrable on $M$, 
 by  lemma 3.1
$\|\nabla f\|^2f^{\alpha}$, and so  $\|\nabla f\|^2f^{\alpha-1}$,
are integrable, and under the condition of integrability  of
$\frac{\sin\theta}{\cos^{\alpha +2}\theta}\|\bar R\|$ we obtain
the integrability of $\Delta f \, f^{\alpha}$. Making $s\to +\infty$,
 from the above inequality we have
$$\int_{M}{(\alpha-1)}
f^{\alpha-2}\|\nabla f\|^2dV=\lim_{s\to +\infty}
\int_{{D_f^-(s)}}{(\alpha-1)}
f^{\alpha-2}\|\nabla f\|^2dV=0$$
and  $f$ is constant, what is impossible
unless $M$ is compact. \qed
\section{Some calibrations}
\renewcommand{\thesection}{\arabic{section}}
\renewcommand{\theequation}{\thesection.\arabic{equation}}
\setcounter{equation}{0}
\subsection{The volume calibration}
We consider in  a Riemannain product $\oM=M\times N$
of two Riemannian manifolds $(M,g_1)$ and $(N,h)$ 
 the volume calibration (1.2), 
and $M$ 
a graph submanifold $F=\Gamma_{\!f}:M\ra M\times N$ of a map
$f:M\ra N$. The graph metric on $M$ is the induced metric 
$g=g_1+f^*h$ by the graph map $\Gamma_{\!f}(p)=(p,f(p))$.
We take $a_i$ a diagonalizing $g_1$-orthonormal basis of
$f^*h$ with eigenvalues $\lambda_1^2\geq \lambda_2^2
\ldots \geq \lambda_m^2\geq 0$. Let $k$ such that $\lambda_k^2>0$
and $\lambda_{k+1}^2=0$, 
and consider the orthonormal system of $T_{f(p)}N$,
 $a_{1+m}, \ldots a_{1+k}$ defined by $df(a_i)=\lambda_i a_{i+m}$,
and extend to an orthonormal basis $ a_{1+m},\ldots, a_{n+m}$.
Then for $i=1,\ldots,m$, $\al=1,\ldots, n$
( where $\lambda_{\al}=0$ for $\al \geq k+1$)
$$X_i=\frac{d\Gamma_{\!f}(a_i)}{
\sqrt{1+\lambda_i^2}}=\frac{a_i+\lambda_i a_{i+m}}{\sqrt{1+\lambda_i^2}},
~~~~~~
X_{m+\al}=\frac{\lambda_{\al} a_{\al}- a_{\al+m}}
{\sqrt{1+\lambda_{\al}^2}}$$
define respectively an orthonormal basis of $(T_{(p,f(p))}\Gamma_{\!f}
,g)$ and of 
$(NM_p, \bar{g})$. The sign of $\lambda_i$ can be chosen such that
$X_i$ is a direct basis of $\Gamma_f$.
Then considering $\Phi$ as a morphism from $T_{(p,f(p))}\Gamma_{\!f}$ 
to $NM_p$ we have 
$$\Phi(X_i)=\cos\theta \lambda_i\,  X_{i+m},~~~~~~~~\cos\theta=
({\Pi_{j}(1+\lambda_j^2)})^{-1/2}$$
and  as a morphism from $TM$,
$\Phi(a_i)
=\cos\theta((df^t df(a_i),-df(a_i))$, where $df^t$ is the adjoint map.
For $B=\sum_{ija}h^a_{ij}X_{a}$,
$a=m+1,\ldots, m+ n$ we have
$$\|B\|^2\geq \sm{\sum}_{ijk}(h^{m+j}_{ik})^2=\sm{\sum}_{i<j,k}
[(h^{m+j}_{ik})^2
+(h^{m+i}_{jk})^2]+\sm{\sum}_{ik}(h^{m+i}_{ik})^2.
$$
If $n=1$, $Q_{\Omega}(B)\geq \|B\|^2$. For $n\geq 2$
 if  we assume $|\lambda_i\lambda_j|\leq 1-\delta$ for $i\neq j$,
 where $0 < \delta \leq 1$,
the quadratic form $Q_{\Omega}(B)$ is also $\delta$-positive. Indeed, 
we have
( \cite{wa2})
\begin{eqnarray*}
Q_{\Omega}(B) &=& \|B\|^2+\sm{\sum}_{ik}\lambda_i^2(h_{ik}^{m+i})^2+
2\sm{\sum}_{k, i<j}\lambda_i\lambda_j h_{jk}^{m+i}h_{ik}^{m+j}\\
&= &\delta \|B\|^2+(1-\delta)\|B\|^2+\sm{\sum}_{ik}\lambda_i^2(h_{ik}^{m+i})^2
+2\sm{\sum}_{k, i<j}\lambda_i\lambda_j h_{jk}^{m+i}h_{ik}^{m+j}\\
&\geq&\delta \|B\|^2+ (1-\delta)\La{(}\sm{\sum}_{i<j,k}[(h^{m+j}_{ik})^2
+(h^{m+i}_{jk})^2]+\sm{\sum}_{ik}(h^{m+i}_{ik})^2\La{)}\\
&& -(1-\delta)\sm{\sum}_{ik}(h^{m+i}_{ik})^2-2(1-\delta)\sm{\sum}_{i<j,k}
|h^{m+j}_{ik}|\, |h^{m+i}_{jk}|
~\geq~\delta \|B\|^2
\end{eqnarray*}
Moreover, $\cos^2\theta\, (\sm{\sum}_i \lambda_i^2) \leq \sin^2\theta
\leq (m-1)\cos^2\theta\, (\sm{\sum}_i\lambda_i^2)$, for 
$$\sin^2\theta=1-\cos^2\theta=\cos^2\theta(\Pi_i(1+\lambda_i^2)-1)=
\cos^2\theta(\sm{\sum}_i\lambda_i^2+\sm{\sum}_{i<j}\lambda_i^2\lambda_j^2
+\ldots).$$
We denote by $\lnab{}df$ the Hessian  of
$f:(M,g_1)\ra (N,h)$.
Let $g_{ij}=g(a_i,a_j)=\delta_{ij}(1+\lambda_i^2)$ and
consider the section $W$ of $f^{-1}TN$ and the vector field
$Z_1$ of $M$:
$$
W= trace_g \lnab{}df= \sm{\sum}_{ij}g^{ij}\lnab{}df(a_i,a_j),~~~~~~~
Z_1=\sm{\sum}_{st}g^{st}h(W,df(a_s))a_t. $$
Then (see \cite{[S2]})  $mH=(-Z_1,W-df(Z_1))$.
Assuming $H$ is parallel with $c=\|H\|$, the vector field $Z$
we used in the proof of theorem 1.2 can be expressed as
$Z= -\frac{\cos\theta}{m}Z_1$.
We have the relations $\|Z_1\|_{g_1}\leq mc$,  $div_{g_1}(Z_1)=m^2c^2$, 
$\|Z\|_g\leq |\sin\theta|c$, $
div_g(Z)=-m c^2\cos\theta$.
Integration on $D$ of $div_{g_1}(Z_1)$  gives theorem 1.1. 
\subsection{The foliation calibration}
Assume 
$\pi:\oM\to N$ is a transversally oriented 
$m$-foliation of an oriented Riemannian manifold
$(\oM^{m+n},\bar{g})$ onto a 
set $N$ with
$\oM=\cup_{y\in N}M_y$, where the leaf at $y$, $M_y=\pi^{-1}(y)$, is 
 an oriented $m$-submanifold, and such that
 for each $x\in \oM$, we have a split $T_x\oM=T_x\oM^v\oplus T_x\oM^h$
where the vertical space $T_x\oM^v=T_x(M_{\pi(x)})$
 is orthogonal to the horizontal space $T_x\oM^h$, defining
smooth oriented vector subbundles of $T\oM$.
If $N$ is a smooth $n$-manifold,
$\pi$ is a fibration if $\forall x$, 
${ Kern}\, d\pi(x)=T_x\oM^v$ and
$d\pi(x):T_x\oM^h\to T_{\pi(x)}N$ is an isomorphism,
and it is Riemannian if $N$ has a metric $h$ such that
$d\pi(x)$ is an isometry for any $x$.
We define for $X_i$ vector fields of $\oM$
$$\Omega(X_1,\ldots, X_m)=Vol_{\pi(x)}(X_1^v, \ldots, X_m^v)$$
where $Vol_{\pi(x)}$ is the volume element of the leaf at $p=\pi(x)$.
Let $e_a$, $ a=1,\ldots, m+n$,
 be a local orthonormal frame
of $\oM$ with
$e_i$, $  i=1,\ldots, m$,  vertical and $e_{\alpha}$, $\alpha=m+1,\ldots, m+n$,
horizontal. We denote by $B^v_x(e_j,e_i)=(\bar{\nabla}_{e_j}e_i)^h$
 the second fundamental form of the leaf  $M_{\pi(x)}$ and $H^v_x$
the mean curvature at $x$. We assume $e_i$ is a direct basis of the leaf.
Now we have
$$\bar{\nabla}_{e_a}\Omega(e_{\alpha},e_1,\ldots,\hat{e}_i, \ldots, e_m)
=-\Omega(\bar{\nabla}{e_a}e_{\alpha},e_1, \ldots, \hat{e}_i, \ldots, e_m)
=(-1)^{i}\bar{g}(\bar{\nabla}e_{a}e_{\alpha},e_i).$$
From this equality and other similar ones, follows the following lemma:
\begin{lemma}
All components of $\bar{\nabla}\Omega$ and
of  $d\Omega$ vanish
except for the following where $i,j\leq m$, $\alpha,\beta\geq m+1$ 
$$ \begin{array}{l}
\bar{\nabla}_{e_j}\Omega\, (e_{\alpha},e_1, \ldots, \hat{e}_i, \ldots, e_m)
=(-1)^{i+1}\bar{g}(B^v(e_j, e_i), e_{\alpha}) \\
\bar{\nabla}_{e_{\beta}}\Omega\, ( e_{\alpha},
e_1, \ldots, \hat{e}_i, \ldots, e_m)
=(-1)^{i}\bar{g}(\bar{\nabla}_{e_{\beta}}e_{\alpha}, e_i)\\
d\Omega\, (e_{\alpha}, e_1,\ldots, e_m) =- m\, \bar{g}(H^v, e_{\alpha})\\
d\Omega\, (e_{\alpha}, e_{\beta},e_1,\ldots,\hat{e}_i, \ldots,  e_m) = 
(-1)^i\bar{g}([e_{\alpha},e_{\beta}], e_i)
\end{array}$$
\end{lemma}
\noindent
So, $d\Omega=0$ iff $H^v=0$ and  $[T\oM^h, T\oM^h]\subset T\oM^h$, and
$\bar{\nabla}\Omega=0$ iff $B^v=0$ and
 $(\bar{\nabla}_{T\oM^h}T\oM^h)^v=0$.
Therefore, we can conclude: 
\begin{proposition} $ \Omega$ is a calibration in $\oM$, that is $d\Omega=0$,
if and only if  the leaves are minimal and the horizontal subspace 
defines an integrable
distribution $T\oM^h$ of rank $n$. 
In this case $\Omega$ calibrates the leaves and
consequently they are stable minimal submanifolds.
It defines a parallel calibration if and only if the leaves 
and the integral submanifolds of $T\oM^h$ are all totally geodesic. 
\end{proposition}
\begin{corollary}
If $n=1$ and $\pi$ defines a foliation of $\oM$ by minimal
hypersurfaces,
then $\Omega$  defines
a closed calibration and  the leaves are stable.
$\Omega$ is  a parallel calibration if and only if
the leaves are totally geodesic and 
the unit normal $\bar{X}$ to the leaves satisfies $\bar{\nabla}_{
\bar{X}}\bar{X}=0$. This is the case when for some function 
$f:M\to \mathbb{R}$,
$\bar{Y}=f\bar{X}$ defines a nondegenerated Killing vector field
with $f$ constant along each leaf.
\end{corollary}
\noindent
\em Proof. \em  We only have to prove the last statement.
From integrability of the leaves, one easily sees that
$0=\bar{g}([e_i,e_j],\bar{Y})=2\bar{g}(e_i, \bar{\nabla}_{e_j}\bar{Y})
=-2\bar{g}(B^v(e_i,e_j),\bar Y)$. Moreover, $f\bar{g}(\bar{\nabla}_{\bar{X}}
\bar{X},e_i)=-\bar{g}(\bar{\nabla}_{e_i}\bar{Y},\bar{X})=-df(e_i).$
\qed\\[3mm]
If $n=1$ the mean curvature and the second fundamental form
of the leaves have been studied  by Barbosa, Kenmotsu, Oshikiri,
Bessa and Montenegro in
\cite{[BKO],[BBM]}, for the general case $\Omega$  closed.
The above corollary 5.1, 
using  the fact $\Omega$ is a closed calibration, 
gives an elementary proof of the stability of
the leaves (see introduction) proved in \cite{[BBM]} 
using more classical stability arguments
involving eigenvalue problems and maximum principles. 
\begin{corollary} Let $(N,h)$ be a Riemannian manifold
and $G$  a Lie group and $N'\subset G$ a subset 
that acts transitively and freely on N as a
group of isometries, and
$f:M\to N$  a smooth map defining a minimal graph
$\Gamma_f$ of $M\times N$. If the orbit 
$(Id_M\times N')\, NM$  of the normal bundle 
of $\Gamma_{f}$ defines an integrable distribution of $M\times N$, 
then $\Gamma_f$ is stable. Furthermore, if $G=N=N'$, and $\g=T_{e}G$,
$e$ the identity element,
then
the later condition is equivalent to
$\psi=-df^t:\g\to C^{\infty}(TM)$ is a Lie algebra homomorphism.
In particular all minimal graphs in $M\times \mathbb{R}$ are stable.
\end{corollary}
\noindent
 \em Proof. \em  This argument is used in \cite{[BBM]} for
$N=G=\mathbb{R}$ to construct a foliation in $M\times N$.
Let $af(p)=L_a(f(p))$, where $a\in N'\subset G$ acts on the left of $N$,
and $a$ on $T_yN$ acts  as $(L_a)_{*y}$.
Since the action is free and transitive then $M\times N=\cup_{a\in G}
\Gamma_{af}$ is a foliation, and since each $a$ is an isometry,
$\Gamma_{af}$ is a minimal submanifold.
Now, the normal bundle of $\Gamma_{af}$ is just given by
$(Id \times (L_a)_*)(NM)=\{(-df_p^t(Y),(L_a)_{*f(p)}Y)
: Y\in T_{f(p)}N \}$,  and  we get the integrability
condition on the orbit. If $G=N=N'$, this
means  
$\{(-df_p^t(\tilde{Y}_{f(p)}),\tilde{Y}_{y}), \tilde{Y}\in \g, y\in G, 
p\in M\}$ is an integrable distribution, where $\tilde{Y}_y=(L_y)_{*e}(
\tilde{Y}_e)$.
We easily see this is  equivalent to
 $[\psi(\tilde{X}),\psi(\tilde{Y})]_M=\psi([\tilde{X},\tilde{Y}])$, 
where $\psi(\tilde{X})$
defines a vector field on $M$, that at $p\in M$
values $-df_p^t(\tilde{X}_{f(p)})$.
\qed\\[3mm]
\em Remark. \em
In the previous corollary, if $m=n$,  $M=\mathbb{R}^n$ and
$f=\nabla \phi$ where $\phi:\mathbb{R}^n\to \mathbb{R}$ is a smooth
function, defining a minimal Lagrangian graph in $\mathbb{R}^n
\times \mathbb{R}^n$, then $\psi(\partial_i)=
Hess\, \phi^{\sharp}(\partial_i)$.
The stability condition means $D^3\phi(Hess\,\phi^{\sharp}(\partial_i),
\partial_j,\partial_k)
$ $=D^3\phi(Hess\,\phi^{\sharp}(\partial_j),\partial_i,\partial_k)$. This seems
to be quite different from the condition of special
Lagrangian graphs. Minimal Lagrangian graphs
are calibrated by the special Lagrangian calibration (\cite{[HL]}), 
and  may not
 satisfy the above condition, that is related to a different calibration.
\begin{corollary} Any minimal submanifold  $M^m$ of $\mathbb{R}^{m+n}$, 
with $n=1$,
is locally  stable and locally a calibrated submanifold.
The same also holds for $n\geq 2$ if  for a local representation
as a graph of a map $f:\mathbb{R}^m\to \mathbb{R}^n$,
$\psi=-df^t:\mathbb{R}^n\to C^{\infty}(T\mathbb{R}^m)$ 
is a Lie algebra homomorphism.
\end{corollary}
\noindent
\em Proof. \em 
Let $p\in M$
and $D$ a compact neighbourhood of $p$ where it is defined
$\nu_{\alpha}$,  an orthonormal frame
of $NM$ defined on $D$. We identify $L=T_pM$ with
 $\mathbb{R}^{m}$. Now $d\pi_L(p)$ is the identity map
of $T_pM$ where
$\pi_L(q)=q-\sum_{\alpha}\bar{g}(q,\nu_{\alpha}(p))\nu_{\alpha}(p)$ and
 we may assume  $\pi_L$ is a diffeomorphism.
Thus,  $D$ is the graph of a map $f:\mathbb{R}^m\to\mathbb{R}^{n}$.
\qed\\[3mm]

Next we generalize
the main result of \cite{[BKO]} (proposition 2.14)
for any codimension.
We assume $T\oM^h$ is an integrable distribution and 
 consider the
maximal horizontal integrable $n$-submanifold $\Sigma$ 
passing at a given pont $x\in \oM$, with second
fundamental form  
$B^h_x(e_{\alpha},e_{\beta})=(\bar{\nabla}_{e_{\alpha}}e_{\beta})^v$, 
and mean curvature 
$H^h_x=\frac{1}{n}\sum_{\alpha}B^h_x(e_{\alpha},e_{\alpha})$.
Let $s^v$ be the scalar curvatures of the leaves and $s^h$ the
ones of the horizontal integrable submanifolds,
and  $\overline{Ricci}$ and $\bar{s}$ the Ricci tensor and scalar curvature
of $\oM$, respectively. $H^h$
defines a vertical vector field of $\oM$, and consider
its divergence in a fiber and  in $\oM$ respectively:
$$ {div}(H^h)=\sum_i\bar{g}(\bar{\nabla}_{e_i}H^h, e_i),~~~~
\overline{{div}}(H^h)=\sum_a\bar{g}(\bar{\nabla}_{e_a}H^h, e_a)$$
Similarly for the horizontal vector field $H^v$ we may take
its divergence along
an horizontal integrable submanifold $\Sigma$, 
$div_{\Sigma}(H^v)=\sum_{\alpha}\bar{g}(\bar{\nabla}_{e_{\alpha}}
H^v, e_{\alpha})$.
\begin{lemma} Assuming $T\oM^h$ is an integrable distribution
\begin{eqnarray*}
\overline{{div}}(H^h) &=& {div}(H^h) -n\|H^h\|^2\\[1mm]
\overline{{div}}(H^v) &=& {div}_{\Sigma}(H^v) -m\|H^v\|^2\\[1mm]
n\, {div}(H^h) + m \, div_{\Sigma}(H^v)
&=&
\sm{\sum}_i\overline{Ricci}(e_i,e_i) -s^v +m^2\|H^v\|^2
+\|B^h\|^2\\
&=&\sm{\sum}_{\alpha}\overline{Ricci}(e_{\alpha},e_{\alpha}) -s^h +
n^2\|H^h\|^2 
+ \|B^v\|^2.\\[1mm]
2n\, {div}(H^h) + 2m \, div_{\Sigma}(H^v)&=&\\
\lefteqn{=\bar{s} -s^v-s^h +m^2\|H^v\|^2+n^2\|H^h\|^2
+\|B^h\|^2+ \|B^v\|^2.}~~~~~~~~~~~
\end{eqnarray*}
\end{lemma}
\noindent
{\em Proof. } We may assume at a point $x$, $(\bar{\nabla}_{e_a}e_i(x))^v=
(\bar{\nabla}_{e_a}e_{\alpha}(x))^h=0$, $\forall a,i,\alpha$. The
first two equalities are obtained by an elementary computation. At $x$
\begin{eqnarray*}
n\, div(H^h)&=& \sm{\sum}_{i,\alpha}\bar{g}(
\bar{\nabla}_{e_i}(\bar{\nabla}_{e_{\alpha}}e_{\alpha}
)^v),e_i)=  \sm{\sum}_{i,\alpha}d(\bar{g}((\bar{\nabla}_{e_{\alpha}}e_{\alpha}
)^v, e_i))(e_i)\\
&=& \sm{\sum}_{i,\alpha}d(\bar{g}(\bar{\nabla}_{e_{\alpha}}e_{\alpha}
, e_i))(e_i)=\sm{\sum}_{i,\alpha}\bar{g}(
\bar{\nabla}_{e_i}\bar{\nabla}_{e_{\alpha}}e_{\alpha},e_i)\\
&=& \bar{g}(~\bar{R}(e_{\alpha},e_i)e_{\alpha}
+\bar{\nabla}_{e_{\alpha}}\bar{\nabla}_{e_{i}}e_{\alpha}
+\bar{\nabla}_{[e_i,e_{\alpha}]}e_{\alpha}~, ~ e_i~).
\end{eqnarray*}
and
$$\begin{array}{l}
\sm{\sum}_i\bar{g}(\bar{\nabla}_{e_{\alpha}}\bar{\nabla}_{e_{i}}
e_{\alpha},e_i)
=\sm{\sum}_i-\bar{\nabla}_{e_{\alpha}}(\bar{g}(e_{\alpha},
\bar{\nabla}_{e_{i}}e_i))=
-\bar{g}(e_{\alpha}, m\bar{\nabla}_{e_{\alpha}}H^v)\\
\bar{g}(\bar{\nabla}_{\nabla{e_i}e_{\alpha}}e_{\alpha}~, ~ e_i~)=
\sm{\sum}_j \bar{g}(\bar{\nabla}_{{e_j}}e_{\alpha}~, ~ e_i~)
\bar{g}(\bar{\nabla}_{{e_i}}e_{\alpha}~, ~ e_j~)\\
\bar{g}(\bar{\nabla}_{\nabla{e_{\alpha}}e_{i}}e_{\alpha}~, ~ e_i~)=
\sm{\sum}_{\beta} \bar{g}(\bar{\nabla}_{{e_{\beta}}}e_{\alpha}~, ~ e_i~)
\bar{g}(\bar{\nabla}_{{e_{\alpha}}}e_{i}~, ~ e_{\beta}~)
\end{array}
$$
leading to
$
\sm{\sum}_{i,\alpha}
\bar{g}(\bar{\nabla}_{[e_i,e_{\alpha}]}e_{\alpha}~, ~ e_i~)
=\|B^h\|^2+\|B^v\|^2$, 
and consequently, 
$$n\, div(H^h)=\sm{\sum}_{\alpha, i}\bar{R}(e_{\alpha},e_i,e_{\alpha},e_i)
+\|B^h\|^2+\|B^v\|^2-m\,div_{\Sigma}(H^v). $$
Now 
$\sum_{i\alpha}\bar{R}(e_{\alpha},e_i,e_{\alpha},e_i)
=\sum_i\overline{Ricci}(e_i,e_i)-\sum_{ij}\bar{R}(e_j,e_i,e_j,e_i)$
and using Gauss equation with respect to a leaf, we obtain the first expression
for $n\, div(H^h)$. Writing
$\sum_{i\alpha}\bar{R}(e_{\alpha},e_i,e_{\alpha},e_i)
=\sum_{\alpha}\overline{Ricci}(e_{\alpha},e_{\alpha})
-\sum_{\alpha\beta}\bar{R}(e_{\alpha},e_{\beta},e_{\alpha},e_{\beta})$
and applying Gauss equations with respect to $\Sigma$ we get the
second expression for $n\, div(H^h)$.
Summing the previous two expressions we obtain the last one.
\qed\\[2mm]
\begin{proposition}
Assume $\oM$ is closed,
 $T\oM^h$ is an integrable distribution and $H^v$ is a divergence free
vector field along each horizontal integral submanifold.
 Then $H^v=0$, i.e.\ the leaves are minimal and $d\Omega=0$.
Furthermore:\\[1mm]
$(1)$ If $H^h$ is also a divergence free
vector field along each leaf, the horizontal integral submanifolds
are  minimal as well.\\[1mm]
$(2)$ If $\overline{Ricci}\geq 0$ and
  $s^v\leq 0$ for all the leaves, 
  the horizontal integrable submanifolds
 are totally geodesic and $ \overline{Ricci}$ vanish in the
direction of all leaves and $s^v=0$.\\[1mm]
$(3)$  If $\overline{Ricci}\geq 0$ and $s^h\leq 0$ for 
all horizontal integral submanifolds $\Sigma$, then they
 are minimal and the leaves of $\pi$ are
totally geodesic and $ \overline{Ricci}$ vanish in the
direction of all horizontal vector fields and $s^h=0$.\\[1mm]
$(4)$ If $\bar{s}\geq 0$ and $s^v+s^h\leq 0$,
all the leaves and the horizontal integrable submanifolds
 are totally geodesic, and $\bar{s}= s^v+s^h=0$.\\[2mm]
Consequently,  if $\bar{s}\geq 0$ ( $\overline{Ricci}\geq 0$ resp.)
and $\bar{s}>0$ ($\overline{Ricci}>0$ resp.)
 at some point $x$,
then either the leaf $M_{\pi(x)}$ or (and resp.) the
horizontal integral submanifold at some  $y\in \pi^{-1}(x)$
must have positive scalar curvature somewhere. 
In particular, in the later case, $n\geq 2$.
\end{proposition}
\noindent
\em Proof. \em 
If $div_{\Sigma}(H^v)=0$ for each $\Sigma$,
 integration of $\overline{div}(H^v)$ on $\oM$  of the second
formula of previous lemma gives  $H^v=0$. Similarly for (1), using the
first formula. 
To prove (2) (3) and (4) we   integrate the last three formulas of
lemma 5.2, with $H^v=0$,  along each leaf $M_{\pi(x)}$, 
that is compact.\qed\\[3mm]
\noindent
If $n=1$ and $H^v$ is constant, with the same constant for all
fibers, then $H^v$ is a divergence free vector along each horizontal
integral submanifold, giving the case in  \cite{[BKO]}.

\subsection{The K\"{a}hler calibration} On a K\"{a}hler  manifold
$(\oM,J,\bar{g})$ with K\"{a}hler form $w(X,Y)=g(JX,Y)$ it is defined
the K\"{a}hler  calibration $\Omega=\frac{w^k}{k!}$, that calibrates the
complex submanifolds of complex dimension $k$.
If 
 $\theta_1,\ldots,\theta_k$
are the K\"{a}hler angles of $M$, $\cos\theta_1\geq \ldots\geq \cos\theta_k\geq 0$ 
and $e_a=X_i,Y_i$ a diagonalizing
o.n basis of $F^*w$, that is $F^*w(X_i,X_j)=F^*w(Y_i,Y_j)=0$,
$F^*w(X_i,Y_j)=\cos\theta_i \delta_{ij}$, then
$$\cos\theta=\epsilon \cos\theta_1\ldots\cos\theta_k,~~~~~\epsilon=\pm 1$$
and $\Phi(X_i)=-\epsilon \cos\theta(J(\frac{Y_i}{\cos\theta_i}))^{\bot}$, ~
$\Phi(Y_i)=\epsilon \cos\theta(J(\frac{X_i}{\cos\theta_i}))^{\bot}$. 
A submanifold $M$ is said to have   equal K\"{a}hler angles,
if $\epsilon=1$ and $\theta_i=\vartheta$ $\forall i$ (see \cite{[SValli]}). 
It is a complex (resp.\ Lagrangian) submanifold iff
$\cos\vartheta=1$ ( resp. $\cos\vartheta=0$).
We assume $M$ and $\oM$ are  of real and complex dimension 4 respectively
and $M$ has equal K\"{a}hler angles. We recall that
 $(F^*w)^{\sharp}:TM\ra TM$ and 
$\Phi:TM\ra NM$ (with respect to
the K\"{a}hler calibration)  are conformal morphisms with
coefficient of conformality $\cos^2\vartheta$ and 
$\sin^2\vartheta\cos^2\vartheta$, respectively.
Note that  $\Phi=
- \Phi'\circ (F^*w)^{\sharp}$ with $\Phi'(X)=(JX)^{\bot}$, given in
\cite{[SValli],[S4]}. We have $\|\Phi'(X)\|^2=\sin^2\vartheta\|X\|^2$.
We can write $(F^*w)^{\sharp}=\cos\vartheta J_{w}$
where  $J_w$ is the almost complex structure
of $M$, defined where $\cos\vartheta\neq 0$
by $J_w(X_i)=Y_i$. 
Similarly we get a polar decomposition
for $w^{\bot}=\cos\vartheta J^{\bot}$, the restriction of $w$ to the 
normal bundle.
The orthonormal frame of the normal bundle
$U_i=\Phi'(\frac{Y_i}{\sin\vartheta})$, $V_i=J^{\bot}U_i
=\Phi'(\frac{X_i}{\sin\vartheta})$  diagonalizes $w^{\bot}$.
We have
$$\begin{array}{ccl}
\Omega(X_{i},Y_j,U_k,V_s) &=& (\cos^2\vartheta+\sin^2\vartheta
\delta_{ik})\delta_{ij}
\delta_{ks}\\
\Omega(X_1,X_2,V_1,V_2) & =&\Omega(Y_1,Y_2,U_1,U_2)=-\sin^2\vartheta
\end{array}$$
and all the other components of $\Omega$ in this basis vanish. Then it follows
the condition $\tilde{Q}_{\Omega}(B)\geq \delta\|B\|^2$ 
is very restrictive.
 Consider the complex and anticomplex parts of $B$
with respect the almost complex structures $J_w$ of $TM$ and
$J^{\bot}$ of $NM$:
$B^c(X,Y)=\sm{\frac{1}{2}}(B(X,Y)-J^{\bot}B(J_wX,Y))$, 
$B^a(X,Y)=\sm{\frac{1}{2}}(B(X,Y)+J^{\bot}B(J_wX,Y)). $
Then $\tilde{Q}(B)=\|B\|^2+(\cos\theta)^{-1}(\|B^a\|^2-\|B^c\|^2)+\rho$,
where 
$$|\rho|\leq 4 \sm{\frac{\sin^2\vartheta}{\cos\theta}}\, \sm{\sum}_{a<b,c
}\|B(e_a,e_c)\|\,\|B(e_b,e_c)\|\leq  12\, 
\sm{\frac{\sin^2\vartheta}{\cos\theta}}
\,\|B\|^2.$$
Using this upper bound, for $\cos\theta\in (\sm{\frac{11}{13},
\frac{11}{12}}]$, if 
\begin{equation}
\|B^a\|^2\geq \sm{\frac{ (13-\cos\theta(13-\delta))}{(-11+\cos\theta
(13-\delta))}}\|B^c\|^2,~~~~\mbox{with~~}
0\leq \delta< \sm{\frac{13\cos\theta-11}{\cos\theta}}\leq 1
\end{equation}
we have  $\tilde{Q}(B)\geq \delta\|B\|^2$.
Note that if
$M$ is a complex submanifold,
 $J_w=J^{\bot}=J$, $\sin\vartheta=0$, $B$ is a complex
bilinear form, and $\tilde{Q}(B)=0$. So, calibrated
submanifolds may not be totally geodesic. Theorem 1.5(A) gives (3)
of next proposition
\begin{proposition} 
Let $F:M^{2k}\to \oM^{2k}$ be a  $2k$-submanifold
immersed with parallel mean curvature $H$ and with equal K\"{a}hler angles
into a K\"{a}hler manifold of complex dimension $2k$
and scalar curvature $\bar{s}$. \\
(1)~\cite{[SValli]} Assume $H=0$ and $\oM$ is  Einstein.
 If $k=2$ and $\bar{s}\neq 0$,
then $F$ is either a complex or a Lagrangian  submanifold.
If $k\geq 3$, $\bar{s}<0$, and  $M$ closed,  then $F$ is either complex
or Lagrangian. If $k\geq 3$, $\bar{s}=0$, and $M$ closed,
 then $\theta$ is constant.\\[1mm]
(2)~\cite{[SPort]} If $k=2$, $\bar{s}<0$, $M$ closed and $\|H\|^2\geq
-(\bar{s}/8)\sin^2\vartheta$, then $F$ is either
a complex or a Lagrangian submanifold.\\[1mm]
(3)~ If $k=2$, $\cos\theta>0$, (1.7) holds and (5.1) is satisfied for
some $\delta$ and $\cos\theta\in (\sm{\frac{11}{13},\frac{11}{12}}]$, 
and $M$ is closed, then $F$ is totally geodesic.
\end{proposition}
\noindent
If $\oM$ is a complex space form of sectional
holomorphic curvature $\nu$
then for $U\in NM$,
$\bar{R}(e_a,e_b,e_a,U)=\frac{3\nu}{4} w(e_a,e_b)\bar{g}(\Phi'(e_a),U)$,
 and  $\sum_{ab}\bar{R}(e_a,e_b,e_a,\Phi'(e_b))$
$=\frac{3k\nu}{2}\cos^2\vartheta\sin^2\vartheta$. Note
that $\nu$ has the same sign has $\bar{s}$, but for $k=2$ and
$\nu>0$  (5.1) in (3) does not hold  because
of (1).  Hence, $\delta$-positiveness of $\tilde{Q}$ 
can be expected only when $\bar{s}=0$.
(2) is related to a result obtained by Kenmotsu and Zhou in \cite{KenmZh},
and Hirakawa in \cite{[Hira]}
where a classification of surfaces with parallel mean curvature
in a complex space forms is obtained using the K\"{a}hler angle.
\subsection{ The Quaternionic calibration}
 This calibration is not so well
understood in the literature so we will describe in some detail.
Let $(V,I,J,K,g)$ be an hyper-Hermitean
vector space of dimension $4n$, where $I,J$ are two
anti-commuting $g$-orthogonal structures. For each $x=(a,b,c)\in \mathbb{S}^2$
it is defined a $g$-orthogonal structure
$J_x=aI+bJ+cK$ and its K\"{a}hler form $w_x(X,Y)=g(J_xX,Y)$.
Then $V$ is a right-quaternionic vector space with
$X\zeta= \zeta_0X-J_xX=:J_{\bar{\zeta}}X$
where $\zeta=(\zeta_0,x)\in \mathbb{H}$ and
$J_x$ is extended linearly for $x\in \mathbb{R}^3$.
The right-quaternionic linear group
of isometries  of $V$ is $Sp(n)=Sp(V)=\{\xi\in O(V):
\xi J_x=J_x\xi, \forall x\in \mathbb{S}^2\}\subset SO(V)$.
Let $Sp(1)=\{\zeta \in  \mathbb{H}: |\zeta|=1\}$.
The inclusion
$Sp(V)\cdot Sp(1)=Sp(V)\times Sp(1)/\pm(Id,1)
\subset SO(V)$ is given by
$(\xi,{\zeta})X=\xi(X)\zeta^{-1}.$
Moreover, for $P=(\xi,{\zeta})$
$$
P(J_xX)=\xi(J_xX){\zeta}^{-1}=-\xi(X)\zeta^{-1}\zeta x\zeta^{-1}
=J_{\tau_{{\zeta}(x)}}(\xi(X)\zeta^{-1})=J_{\tau_{{\zeta}(x)}}P(X)
$$
 where $\tau:S^3
\subset \mathbb{H}\ra SO(3)\subset SO(4)$ is the double covering map
$\tau_{\zeta}(v)=\zeta v \bar{\zeta}$.
A subspace $T$
is a \em complex subspace \em if $J_xT\subset T$ for some $x\in \mathbb{S}^2$.
It is a quaternionic subspace if it is $J_x$-complex
$\forall x$.
The  fundamental 4-form of $V$ is defined by
\begin{equation}
\Omega=\sm{\frac{1}{6}}(w_I\wedge w_I
+w_J\wedge w_J+ w_K\wedge w_K).
\end{equation}
For each $X\in V$, let $H_X^0=span\{IX,JX,KX\}$,
$H_X=\mathbb{R} X\oplus H_X^0$.
Each $P\in SO(V)$ acts on $\Omega$ as
$P\Omega(X,Y,Z,W)=\Omega(P^{-1}X,P^{-1}Y,P^{-1}Z,P^{-1}W)$,
and we have
\begin{lemma}
 $H_{\Omega}:= \{P\in SO(V): P\cdot \Omega=\Omega\}
=Sp(V)\cdot Sp(1)$.
\end{lemma}
\noindent
\em Proof. \em
If $P=(\zeta,\xi)\in \mathbb{S}^2 \times Sp(V)$ then
$P=(\xi,\zeta)$ satisfies $P(H_X)=H_{P(X)}$.
Note that $\forall X,Y\in V$,
$w_I(IX,IY)=w_I(X,Y)$, $w_J(IX,IY)=-w_J(X,Y)$.
If $P=(\xi, q)\in Sp(V)\cdot Sp(1)$  one can prove directly
that $P\cdot \Omega=\Omega$ ( see \cite{[Kraines]}).
Then $Sp(V)\cdot Sp(1)\subset
H_{\Omega}$. Now if $P\in H_{\Omega}$,  from the above considerations
$P(H_X)=H_{P(X)}$. Thus,  $\forall x\in \mathbb{R}^3$
$P(J_xX)=J_{A(X,x)}P(X)$,
 with $A(X,\cdot)\in SO(3)$ necessarily in case $\|X\|=1$.
We extend  $A(X,\lambda)=\lambda$, for $\lambda \in
\R\subset \mathbb{H}$.
Since $\forall \lambda\in \R$,
$P(J_x\lambda X)=\lambda P(J_xX)$
we get $A(\lambda X,x)=A(X,x)$. Now  we assume $\|X\|=1$.
From $P(J_{xy}X)=P(J_x(J_yX))$ we have
$A(X,xy)=A(J_yX,x)A(X,y)$. Let $x=\mu x'+\lambda y$, where $x'\bot y$ is a
unit of $\R^{3}$.
Then $A(X,x'\times y)=A(X,x')\times A(X,y)$, and so
$A(J_yX, x)$ $A(X,y)=A(X,xy)=\mu A(X,x'\times y)-\lambda Id=
\mu A(X,x')\times A(X,y)-\lambda Id$, 
implying
$A(J_yX,x)$ $=\mu A(X,x')+\lambda A(X,y)=A(X, \mu x'+\lambda y)=A(X,x).$
Finally let $X,Y$ units with $H_X\oplus H_Y$ and $Z=\frac{X+Y}{\|X+Y\|}$.
Then $H_{P(X)}\oplus H_{P(Y)}$. 
From $P(J_x(X+Y))=J_{A(Z,x)}P(X+Y)$
we get $J_{A(X,x)}P(X)+J_{A(Y,x)}P(Y)=J_{A(Z,x)}P(X)
+J_{A(Z,x)}P(Y)$, and so
$A(X,x)=A(Z,x)=A(Y,x)$.
Then  $A(X,x)=A(x)$ $\forall X$, that
 is $A$ does not depend on $X$. We have proved that
$P(J_xX)=J_{\tau_{\zeta}(x)}P(X)$, where $A=\tau_{\zeta}$
for some $\zeta\in Sp(1)$ (unique up to a sign).
Define $\xi:V\ra V$ by
$\xi(X)=P(X){\zeta}=J_{\bar{\zeta}}P(X)$. Then
$P=(\xi,\zeta)$ and  $\xi\in Sp(V)$.
\qed\\[4mm]
The  fundamental 4-form induces a 
symmetric endomorphism
$\Omega^{\Delta}:\wedge^2V\ra \wedge^2V$, defined by
$\langle\Omega^{\Delta}(X\wedge Y),Z\wedge W\rangle=\Omega(X,Y,Z,W).$
For each oriented 
orthonormal system $B=\{X_1,X_2,X_3,X_4\}$ 
 of $V$, we define the bivectors
 $\Lambda^{\pm}_r=\Lambda^{\pm}_r(B)$, by
{\footnotesize
$$\Lambda^{\pm}_1=\sm{\frac{1}{\sqrt{2}}}(X_1\wedge X_2\pm X_3\wedge X_4),~~
\Lambda^{\pm}_2=\sm{\frac{1}{\sqrt{2}}}(X_1\wedge X_3\mp X_2\wedge X_4),~~
\Lambda^{\pm}_3=\sm{\frac{1}{\sqrt{2}}}(X_1\wedge X_4\pm X_2\wedge X_3).
\\$$}
\noindent
If $X$ is a unit,
$\Lambda_r^{\pm}(X)$ is defined as above w.r.t.
$X_1=X, X_2=IX, X_3=JX, X_4=KX$.
Note that $\Lambda_r^{\pm}(X)=\Lambda_r^{\pm}(J_xX)$ for any 
$x\in \mathbb{S}^2$.
Set for any $X,Y$ and $i=0,1,2,3$, $r=1,2,3$, and 
$\epsilon^r_0=\epsilon_1^3=\epsilon_2^1=\epsilon_3^2=+1$, $
\epsilon_1^1=\epsilon_1^2=\epsilon_2^2=\epsilon_2^3=\epsilon_3^1=
\epsilon_3^3=-1$,
{\small $$
\Theta_i(X,Y)=\sm{\frac{1}{2}}(X\wedge Y+ \epsilon_i^1 IX\wedge IY+
\epsilon_i^2JX\wedge J Y+\epsilon_i^3KX\wedge KY).
$$}
\noindent
satisfying $\sm{\Theta_{s'}(J_xX,J_xY)}\in 
span\{\sm{\Theta_s(X,Y)},s=0,\ldots 3\}$ and
for $J_x=I,J,K$, $\sm{\Theta_{s'}(X,J_xX)}$  either is zero or gives
$\Lambda^{\pm}_s(X)$ for some $s$.
If $H_X,H_Y, H_Z$ are orthogonal quaternionic lines, and
$X,Y$ are units we have for any $x,y,z\in \mathbb{S}^2$
 \begin{equation}
\begin{array}{l}
\sm{\Omega(X,J_xX,J_yX,J_zX)=\langle x,y\times z\rangle,}\\
\sm{\Omega(X,J_xX,Y,J_yY)=\frac{1}{3}\langle x,y\rangle,}\\
\sm{\Omega(X,J_xX,J_yX,Y)=\Omega(X,J_xX,Y,Z)=0}.\end{array}\\
\end{equation}
\noindent
We take an orthonormal basis $X_i$ of $V$ the form
$\{e_{\al},Ie_{\al},Je_{\al},Ke_{\al}\}$, $\al=1,\ldots, n$.
 We have $\Omega^{\Delta}(\xi)=\sum_{i<j}\Omega(\xi, X_i\wedge X_j)
X_i\wedge X_j$.
An orthonormal basis of eigenvectors of $\Omega^{\Delta}$ is given 
by the $2n(4n-1)$
vectors, where $\al,\be=1,\ldots,n$, $r=1,2,3$,
 \begin{equation}\begin{array}{l}
\ft{(1/\sqrt{n})\sum_{\al}\Lambda^+_r(e_{\al}),~~~
\Lambda^-_r(e_{\al}),~~~
1/\sqrt{2}(\Lambda^+_r(e_{n})-\Lambda^+_r(e_{\al})) ~~\al<n}\\
\ft{\Theta_s(e_{\al},Ie_{\be}),~~~
\Theta_s(e_{\al},Je_{\be}),~~~\Theta_s(e_{\al},Ke_{\be})~~~\al<\be,~ 
s=0,1,2,3}
\end{array}\end{equation} 
\noindent
The corresponding eigenvalues,
that range $\{\frac{2n+1}{3},\pm 1, \pm 1/3\}$, are given
as follows where $i=1,2,3$,  $L=id,I,J,K$
{\footnotesize $$\begin{array}{ll}
\Omega^{\Delta}(\Lambda^+_r(e_{\al}))=\Lambda^+_r(e_{\al})
+\sum_{\be\neq \al}\frac{2}{3}\Lambda^+_r(e_{\be}) &
\Omega^{\Delta}(\sum_{\al}\Lambda^+_r(e_{\al}))
=\frac{2n+1}{3}(\sum_{\al}\Lambda^+_r(e_{\al}))\\
\Omega^{\Delta}(\Lambda^+_r(e_{\al})-\Lambda^+_r(e_{\be}))
=\frac{1}{3}(\Lambda^+_r(e_{\al})-\Lambda^+_r(e_{\be}))&
\Omega^{\Delta}(\Lambda^-_r(e_{\al}))=-\Lambda^-_r(e_{\al})\\
\Omega^{\Delta}(\Theta_0(e_{\al},Le_{\be}))=\Theta_0(e_{\al},Le_{\be})&
\Omega^{\Delta}(\Theta_i(e_{\al},Le_{\be}))=-\frac{1}{3}
\Theta_i(e_{\al},Le_{\be})
\end{array}$$}
Given a $k$-dimensional $T$ subspace of $V$  
we consider the restriction
$\Omega^{\Delta}_{T}:\wedge^2T\ra \wedge^2T,$
symmetric endomorphism
with eigenvalues $\al_1,\ldots,\al_{\frac{k(k-1)}{2}}$ that we
call the nonnormalized quaternionic angles of $T$.\\

From now on we restrict our attention when
 $T$ is an oriented four dimensional subspace with direct orthonormal basis
$(X_1,X_2,X_3, X_4)$ and $V$ eight dimensional.
\begin{proposition} The fundamental form 
$\Omega$ defines a calibration that
calibrates the quaternionic 4 dimensional subspaces.
The quaternionic angle of an oriented
$4$-dimensional subspace $T^4$ is defined by the number
$\cos\theta=\Omega(X_1,X_2,X_3,X_4)\in [-1,1].$
$T$ and $T^{\bot}$ have the same quaternionic angle and  there are only
two eigenvalues $\alpha_i=\pm\cos\theta$ each with multiplicity three.
\end{proposition}
\noindent
\em Proof. \em 
To see that $\Omega$ is a calibration, 
we set $\phi(x)=\langle
(w_x)_{|T}\wedge (w_x)_{|T}, Vol_T\rangle$, where $(w_x)_T$ is
 the restriction of $w_x$ to $T\times T$. If $\cos\theta^x_{1}, $  
$\cos\theta^x_2$,
with  $\theta_i^x\in [0,\frac{\pi}{2}]$, are the
$J_x$-K\"{a}hler angles of $T$ w.r.t. $J_x$, then
for any o.n.b. $x,y,z$ of $\R^3$,
\begin{eqnarray}
\cos\theta=\Omega(T) &=&\sm{\frac{1}{6}}(\phi(x)+\phi(y)+\phi(z))\nonumber\\
&=&
\sm{\frac{1}{3}}(\epsilon_x\cos\theta^x_1\cos\theta^x_2+
\epsilon_y\cos\theta^y_1\cos\theta^y_2+
\epsilon_z\cos\theta^z_1\cos\theta^z_2)~~~~
\end{eqnarray}
where $\epsilon_u=\pm 1$ depending if $(w_u)_T$ defines the same
or the opposite orientation of $T$.
From (5.5) we see that $|\phi(u)|\leq 2$ and so
$|\Omega(T)|\leq 1$, and $\Omega(T)= \pm 1$ iff $\epsilon_
{u}=\pm 1$  and $\cos\theta^u_{s}=1$, $s=1,2$,
that is $T$ is a $J_u$-complex subspace,  $\forall u=x,y,z$,
or equivalently, $T$ is a quaternionic
subspace. 
Since the $J_s$-K\"{a}hler angles of $T$ and the ones of the orthogonal
complement of $T^{\bot}$
are the same, then $\Omega(T)=\Omega(T^{\bot})$.\qed
\begin{proposition} If $V$ is 8-dimensional and 
$T$ is a 4-dimensional subspace $J_x$-complex for some $x$ then
$\frac{1}{3}\leq \cos\theta(T)\leq 1$, with equality to
$\frac{1}{3}$ if and only if $T$ is a totally complex subspace,
that is $T$ is a $J_y$-Lagrangian subspace $\forall y\bot x$.
Moreover, if two complex subspaces $T$ and $T'$ of $V$ have the same
quaternionic angle, then there exist an element $P\in Sp(V)\cdot Sp(1)$
such that $T'=P(T)$.
\end{proposition}
\noindent
\em Proof. \em  
If $T$ is a $J_x$-complex subspace of $V$ then
$T$ is Cayley subspace of $(V,g,J_y)$ for all $y\in \mathbb{S}^2$. To see this
we  take an orthonormal basis of $T$ of the form 
$B=\{X_i\}=\{X,J_xX,Z, J_xZ\}$. We have
$g( J_yJ_xX,X) =-\langle y,x\rangle =g( J_yJ_xZ,Z)$,
 $g(J_yJ_xZ,X)
=-g( J_yX, J_xZ)$,
$ g(J_yZ, J_xX)=-g( J_{y\times x}X, Z ),$
and  $g( J_yJ_xZ,J_xX)
=g(J_yX, Z)$.
A basis for the self-dual 2-forms on $T$ is given by $
J^B_r=\Lambda_r^+(B)$.
Then we see that $ (w_y)_{|T}=\cos\theta^y(p)J_v^B$
where
$ v=\frac{1}{t}(<y,x>,g( J_yX,Z),g( J_{y\times x}X,Z))\in \mathbb{S}^2$
 and
$\cos\theta^y(p)=t=\|(J_yX)^{\top}\|$,
proving that $(w_y)_{T}$ is self-dual, that is $*(w_y)_{|T}=(w_y)_{|T}$.
This is just the same as to say the $J_y$-K\"{a}hler angles
of $T$ are equal, that is $T$ is a Cayley subspace.
Therefore $\theta^u_1=\theta^u_2=:\theta^u$ and so
$\cos\theta=\frac{1}{3}(1 +\cos^2\theta^y+\cos^2\theta^z)\geq \frac{1}{3}$,
with equality if and only if $\cos^2\theta^y=\cos^2\theta^z=0$,
that is $T$ is a $J_u$-Lagrangian subspace for any $u\bot x$.
If $T$ and $T'$ have the same quaternionic angle, we use the
canonical frames given in (5.8) below for $T$ and for $T'$
and define $P$ by $P(B)=B'$, $P(B^{\bot})={B'}^{\bot}$.
Then  $PJ_u=J_{u'}P$ for $(u,u')=(x,x'),(y,y')$ or $ (z,z')$,
and $P(X)=X'$, $P(Y)=Y'$, and taking $\zeta\in Sp(1)$ such that
$\tau_{\zeta}$ maps $(x,y,z)$ to $(x', y',z')$, we get 
$\xi(\cdot)=P(\cdot)\zeta^{-1}\in Sp(V)$, what proves that
$ P=(\xi,\zeta)\in Sp(V)\cdot Sp(1)$.
\qed\\

Some further algebraic considerations.
Let $T$ be a Euclidean space of dimension 4. For each  linear map
$l:T\ra T$ we define $\wedge^2l:\wedge^2T\ra \wedge^2T$,
$\wedge^2l(u\wedge v)=l(u)\wedge l(v)$.
If $\lambda_i\geq 0$ are the eigenvalues of $\sqrt{l^Tl}$
(also called the  singular  values of $l$)
and $B=\{e_i\}$ a corresponding orthonormal basis of eigenvectors, let
$\tilde{e}_i$ defined by $l(e_i)=\lambda_i\tilde{e}_i$
whenever $\lambda_i\neq 0$, and extend to an orthonormal basis $\tilde{e}_i$
of $T$. Using Newton inequalities we have
$2\sum_{1\leq i<j\leq 4}\lambda_i\lambda_j \leq 3(\lambda_1^2+
\ldots +\lambda_4^2)$ with equality iff $\lambda_i=\lambda_j$ $
\forall i,j$.
Each direct orthonormal  basis $B$ of $T$ and $B^{\bot}$ of $T^{\bot}$
 define respectively a direct orthonormal basis
$\Lambda^{\pm}_{r}=\Lambda^{\pm}_{r}(B)$
of $\wedge^2_{\pm}T$,
and ${\Xi}^{\pm}_r$ of $\wedge^2_{\pm}T^{\bot}$.
We consider the two hyper-Hermitean structures of $T$
(denoted by $J_s^T$, when we choose one)
$J^B_r=\Lambda_r^+$, 
 $\tilde{J}^B_r=\Lambda_r^-$.
We note the following: If $u,v$ is an o.n. system of vectors of $T$
then
\begin{equation}
|\langle u\wedge v, \Lambda_r^{+}\rangle|\leq 1,
\mbox{~~ with~equality~
to~} 1\mbox{~iff~}v=\pm J_r^B(u).
\end{equation}
and similar for  $\Lambda_r^{-}$.
We define
\[
 Q^{\pm}l = -\sm{\frac{1}{3}}\sm{\sum}_rJ_r^T \circ l\circ J_r^T, ~~~~~~~~
H^{\pm}l = \sm{\frac{1}{4} }(l+ 3Q^{\pm}l)
\]
where $\pm$ depends on $J_r^T=J^B_r$ or $\tilde{J}^B_r$. Note that
$H^{\pm}:Skew(T)\ra sp_1(T)$
gives  the orthogonal projection
of $l$ onto a $(J^T_r)$-hyper-complex linear map (does not depend on the 
oriented basis $B$).
 We also have
$\langle l, Q^{\pm}(l)\rangle =\sum_{r=1,2,3}\frac{4}{3}
\langle \wedge^2l(\Lambda_r^{\pm}), \Lambda_r^{\pm}\rangle$
and  that
$|\langle \wedge^2l(\Lambda_1^{\pm}), \Lambda_1^{\pm}\rangle|
\leq \frac{1}{2}(\lambda_1\lambda_2+\lambda_3\lambda_4)$,
$|\langle \wedge^2l(\Lambda_2^{\pm}), \Lambda_2^{\pm}\rangle|
\leq \frac{1}{2}(\lambda_1\lambda_3+\lambda_2\lambda_4)$,
$|\langle \wedge^2l(\Lambda_3^{\pm}), \Lambda_3^{\pm}\rangle|
\leq \frac{1}{2}(\lambda_1\lambda_4+\lambda_2\lambda_3)$,
and $\|l\|^2=\lambda_1^2+\lambda_2^2+\lambda_3^2+\lambda_4^2$.
Moreover,
\begin{eqnarray*}
9\|Q^{\pm}l\|^2
&=& 3\|l\|^2-2\sm{\sum}_r\langle J_r^T \circ l\circ J_r^T,l\rangle
={3}\|l\|^2 +6\langle Q^{\pm}l,l\rangle,
\end{eqnarray*}
and so $0\leq 16\|H^{\pm}l\|^2
=\La{(}\|l\|^2 +6\langle Q^{\pm}l,l\rangle +
9\|Q^{\pm}l\|^2\La{)}={4}\|l\|^2
+12\langle Q^{\pm}l,l\rangle =16 \langle H^{\pm}l,l\rangle$.
Consequently
$3\langle Q^{\pm}l,l\rangle \geq -\|l\|^2$.
If equality  holds,
then $H^{\pm}l=0$, and so $l\in
sp_{1}(T)^{\bot}=\wedge^2_{+}T=span\{J^T_r\}$.
 Newton inequalities and (5.6) prove that
\begin{equation}
\begin{array}{c}
-\frac{1}{3}\|l\|^2\leq \langle Q^{\pm}l,l\rangle
\leq \|l\|^2,\\[1mm]
\langle Q^{\pm}l,l\rangle=\|l\|^2 \mbox{~iff~} l
\mbox{~is~hyper-complex},\\[1mm]
\langle Q^{\pm}l,l\rangle=-\frac{1}{3}\|l\|^2 \mbox{~~iff~~}
l\in \wedge^2_{+}T=span\{J^T_s\}
\end{array}
\end{equation}
Furthermore, if $l$ is hyper-complex then $l$ is conformal.
 The singular values
of $\wedge^2l$ are $\lambda_i\lambda_j$ for $i<j$.
We can split $\wedge^2l=\wedge^+_+l\oplus \wedge^+_-l
\oplus \wedge^-_+l\oplus \wedge^-_-l$, $\wedge^{\pm}_{+}l:\wedge^2_{+}T
\ra \wedge^2_{\pm}T$, $\wedge^{\pm}_{-}l:\wedge^2_{-}T
\ra \wedge^2_{\pm}T$.
$\wedge^2 l$ is self dual (resp. anti-self-dual),
 i.e. $\wedge^2 l *=*\wedge^2 l$ (resp. $\wedge^2 l *=-*\wedge^2 l$)
iff the anti-self dual part $\wedge^-_+l\oplus\wedge^+_-l$ vanish
(resp. the self-dual part $\wedge^+_+l\oplus\wedge^-_-l$ vanish), iff
either $\wedge^2l=0$,  what means at least 3 of the singular values
vanish, or $l$ is an orientation preserving (resp.\ reversing) conformal
 isomorphism.
\\[2mm]

If $(\oM,g,Q)$ is a quaternionic-K\"{a}hler manifold
of real dimension $4n$ and fundamental
form $\Omega$, the quaternionic $4m$-submanifolds,
 are necessarily totally geodesic (\cite{[Gray]}).
Some attention have been drawn to
a more general  type of submanifolds, the almost complex submanifolds
in the quaternionic context, and their minimality have been
studied. This includes the quaternionic submanifolds as
well the totally complex or the K\"{a}hler submanifolds.
See for example \cite{[AM]} and their references, where some examples
can be found. Most of these  submanifolds are also 
proved to be totally geodesic.
We will show some use of the quaternionic angle
in the study of almost complex submanifolds with parallel mean curvature.

An immersed submanifold $F:M\ra \oM$ is an
almost complex
submanifold if there exist a smooth section $J_M:M\ra Q$ such that, for
each $p\in M$,
$J_M(p)(T_pM)\subset T_pM$.
If $n=2$ and $m=1$,
the quaternionic angle satisfies
$\frac{1}{3}\leq \cos\theta\leq 1$ with equality to $\frac{1}{3}$
at totally complex points and to 1 at quaternionic points.
Since $\cos\theta\geq \frac{1}{3}$ we conclude from theorem 1.4:
\begin{proposition} If $(\oM,g,Q)$ is a quaternionic-K\"{a}hler manifold
of real dimension 8 and $M$ is an almost complex complete submanifold
of real dimension 4 and  with parallel mean curvature
and $Ricci^M\geq 0$, then $M$ is a minimal submanifold.
\end{proposition}
For an almost complex  4-dimensional submanifold,
 $\Phi:TM\ra NM$ is a conformal
morphism with coefficient of conformality $(1-\cos\theta)(\cos\theta-
\frac{1}{3})$
 (\cite{[S4]}). To see this we first note that
we can take \em canonical  orthonormal basis \em
$B$ of $T_pM$ and $B^{\bot}$ of $NM_p$ of the form
\begin{equation}
\begin{array}{ccl}
B &=& \{X,J_{x}X, cJ_{y}X+sY, cJ_{z}X+sJ_{x}Y\}=\{X_k\}\\
B^{\bot} &=& \{J_yY,J_{z}Y, cY-sJ_{y}X, cJ_{x}Y-sJ_{z}X\}=\{U_i\}
\end{array}
\end{equation}
where  $c^2+s^2=1$,
 $x,y,z=x\times y$ is an o.n. basis of $\mathbb{R}^3$ with $J_x=J_M(p)$ and
$Y\in H_X^{\bot}$.
Then $\cos\theta=(1-\frac{2}{3}s^2)$,
$s^2=\frac{3}{2}(1-\cos\theta)$,
$c^2=\frac{3}{2}(\cos\theta-\frac{1}{3})$, and using this basis we
see that $\Phi(B)=-\frac{2}{3}scB^{\bot}$.

Next we use the formula of $\Delta\cos\theta$ to obtain some nonexistence
results for almost complex submanifolds, and in particular to give
 a "calibration"-type proof of the above mentioned result of Gray \cite{[Gray]},
for the case
$n=2$ and $m=1$.
We take for basis of $NM_p$, that is reordering $B^{\bot}$, 
$B'^{\bot}=\{U'_i\}$, $U'_1=U_3$, $U'_2=U_4$, $U'_3=U_1$, 
$U'_4=U_2$, 
and consider the corresponding basis
${\Xi'}^{\pm}_t$,  of $\wedge^2NM_p$.
The matrix of $\Psi: \wedge^2TM \ra \wedge^2 NM$, with respect to
the basis $\Lambda_1^+,\Lambda_2^+,\Lambda_3^+,
\Lambda_1^-,\Lambda_2^-,\Lambda_3^-$ of $\wedge^2TM$ and the basis
 ${\Xi'}_1^+,{\Xi'}_2^+,{\Xi'}_3^+,
{\Xi'}_1^-,{\Xi'}_2^-,{\Xi'}_3^-$ of $\wedge^2NM_p$ 
is given by

\begin{equation}
 \Psi=\left[\begin{array}{cccccc}
\scr{\frac{2}{3}(1+s^2)}&\scr{0}&\scr{0}&\scr{0}&\scr{0}&\scr{0}\\[-1mm]
\scr{0}&\scr{\frac{2}{3}c^2}&\scr{0}&\scr{0}&\scr{0}&\scr{0}\\[-1mm]
\scr{0}&\scr{0}&\scr{\frac{2}{3}c^2}&\scr{0}&\scr{0}&\scr{0}\\[-1mm]
\scr{0}&\scr{0}&\scr{0}&\scr{\frac{2}{3}s^2}&\scr{0}&\scr{0}\\[-1mm]
\scr{0}&\scr{0}&\scr{0}&\scr{0}&\scr{\frac{2}{3}s^2}&\scr{0}\\[-1mm]
\scr{0}&\scr{0}&\scr{0}&\scr{0}&\scr{0}&\scr{-\frac{2}{3}s^2}
\end{array}\right]
\end{equation}
Note that $w_M=\sqrt{2}\Lambda_1^+$  and $w_{NM}=
\sqrt{2}{\Xi'}_1^+$ are the respective K\"{a}hler forms.
$\Psi$ applies $\wedge^2_{\pm}TM$ into $\wedge^2_{\pm}NM$
and denoting the corresponding restriction $\Psi_{\pm}:
\wedge^2_{\pm}TM\ra \wedge^2_{\pm}(NM)$, and defining
$\Psi'_+:=\Psi_+-2(1-\cos\theta)\Psi_0$,
where $\Psi_0:\wedge^2TM\to\wedge^2NM$  is the linear morphism
given by $\Psi_0(w_M)=w_{NM}$ and zero on the
orthogonal complement of $\mathbb{R}w_M$,
then $\Psi'_+$
 and $\Psi_-$ are conformal, with
$\|\Psi'_+(\eta)\|^2=(\cos\theta-\sm{\frac{1}{3}})^2\|\eta\|^2,$
$\forall \eta\in
\wedge^2_+TM,$ and
$\|\Psi_-(\eta)\|^2=(1-\cos\theta)^2\|\eta\|^2,$ $\forall \eta\in
\wedge^2_-TM $.
Thus,
if $M$ is  immersed with no totally complex points
the bundles $\wedge^2_+TM$ and $\wedge^2_+MN$ are isomorphic.
If $M$ is immersed with no quaternionic points, then
$\wedge^2_-TM$ and $\wedge^2_-MN$ are isomorphic.
If there are neither quaternionic  nor totally complex
points, then $\Phi:TM \to NM$ is an isomorphism.\\[3mm]
If $X_i$ is a direct o.n. basis of $T_pM$
and  $Y_i\in NM_p$ are any vectors, then
$$\sm{\sum}_{s<j}\Omega(X_1,\ldots, {Y_s}_{\ti{(s)}},\ldots ,
{Y_j}_{\ti{(j)}},\ldots, X_4)=\sum_r\Omega(\Lambda_r^+, \Lambda_r^+(Y))
-\Omega(\Lambda_r^-, \Lambda_r^-(Y))$$
with $*\Lambda_r^-=-\Lambda_r^-$, and where $\Lambda_r^{\pm}(Y)$,
$Y=(Y_1,\ldots, Y_4)$,  are
formally defined in the same way as $\Lambda_r^{\pm}(B)$.
Thus, we
 consider the two components of $\wedge^2 B$, 
 $\wedge^+_+B:\wedge^2_+T_pM\to \wedge^2_+NM_p$
and  $\wedge^-_-B:\wedge^2_-T_pM\to \wedge^2_-NM_p$, and say that
$\wedge^2 B$ is self-dual iff $\wedge^2 B=\wedge^+_+B\oplus \wedge^-_-B$
and it is anti-self-dual if the self dual part vanish.
Therefore, from lemma 3.1
 \begin{eqnarray}
\Delta \cos\theta &=& -\cos\theta\|B\|^2+
2\langle \Psi_+,\wedge^+_+B\rangle +2\langle \Psi_-,\wedge^-_-B\rangle
\nonumber\\
&&+m\langle\nabla^{\bot}H,\Phi\rangle -\sm{\sum}_{ij}\bar{R}(X_i,X_k,X_i,
\Phi(X_k))
\end{eqnarray}
Now we prove the classic result on quaternionic submanifolds
in \cite{[Gray]} reducing it to
a linear algebra problem:
\begin{proposition} If $M$ is a quaternionic submanifold of 
$\oM$ then $M$ is totally
geodesic.
\end{proposition}
\noindent
\em Proof. \em  We identify
$T_pM=H_{X}\equiv H_{Y}=NM_p$, through the canonical basis 
$B=\{X,J_xX,J_yX,J_zX\}=\{X_k\}$ and $
(B')^{\bot}=\{Y,J_xY,J_yY,J_zY\}$, and set
$l_k=B(X_k,\cdot):T_pM\ra NM_p\equiv T_pM$.
 We have $\Phi=0$ and
$\Psi_-=0$, $\Psi_+=\frac{2}{3}Id$.
 Then (5.10) is
$0= \sum_k \langle Q^{\pm}_B(l_k), l_k\rangle -\|l_k\|^2.$
By (5.7) each $l_k$ is hyper-complex, that is
$B(X_k,J_x Z)=J_x(B(X_k, Z))$, $\forall x\in \mathbb{S}^2$.
Consequently, 
$ B(J_xZ,J_xW)=J_x^2B(Z,W)=-B(Z,W)$, for any $x$. But then
$B(Z,W)=-B(J_yZ,J_yW)=B(J_xJ_yZ,J_xJ_yW)=B(J_zZ,J_zW)$, 
for any o.n. basis $x,y,z=x\times y$ of $\mathbb{R}^3$ and so
$B=0$. \qed\\[3mm]
At $p$
consider the canonical  frames $B=\{X_k\}$  and
$B'^{\bot}=\{U'_k\}$, and
 the linear isometry
$L:NM_p\ra  T_pM,$  $L(U'_k)=X_k$.
We define
$l_k= L\circ B(X_k,\cdot):T_pM\ra T_pM$, and
$l'_k=l_k\circ S'$ and $l''_k=l_k\circ S''$, where $S',S''$ are
orientation preserving isometries, by
\[\begin{array}{c}
l'_k(X_1)= l_k(X_1),~~l'_k(X_2)= l_k(X_2),~~l'_k(X_3)= -l_k(X_3),~~
l'_k(X_4)= -l_k(X_4)\\
{l''}_k(X_1)= l_k(X_1),~~{l''}_k(X_2)= -l_k(X_2),
~~{l''}_k(X_3)= -l_k(X_3),~~
{l''}_k(X_4)= l_k(X_4). \end{array}\]
Then,  $\|l'_k\|=\|{l''}_k\|=\|l_k\|=\|B(X_k,\cdot)\|$, and
we have
$~\wedge^2l'_k(\Lambda^{\pm}_1)=\wedge^2l_k(\Lambda^{\pm}_1)$,
$\wedge^2l'_k(\Lambda^{\pm}_2)=-\wedge^2l_k(\Lambda^{\pm}_2)$,
$\wedge^2l'_k(\Lambda^{\pm}_3)=-\wedge^2l_k(\Lambda^{\pm}_3)$,
$\wedge^2{l''}_k(\Lambda^{\pm}_1)=-\wedge^2l_k(\Lambda^{\pm}_1)$,
$\wedge^2{l''}_k(\Lambda^{\pm}_2)=-\wedge^2l_k(\Lambda^{\pm}_2)$,
$\wedge^2{l''}_k(\Lambda^{\pm}_3)=\wedge^2l_k(\Lambda^{\pm}_3)$.
Set
\begin{eqnarray}
D&=& \sm{\sum_{k}}\la{(}\|l_k\|^2-
\langle Q^+l_k,l_k\rangle \la{)}\geq 0\\
A&=& \sm{\sum_k}\la{(}\langle Q^{+}l'_k,l'_k\rangle +
\sm{\frac{1}{3}}\|l_k\|^2\la{)} \geq 0\\
E&=& \sm{\sum_k}\la{(}
\langle Q^{-}{l''}_k,{l''}_k\rangle + \sm{\frac{1}{3}}\|l_k\|^2\la{)}\geq 0
\end{eqnarray}
Note that $\wedge^2B=\sum_k\wedge^2l_k$ and is antiselfdual
iff $\sum_k\wedge^+_-l_k\oplus \wedge^-_+l_k=0$.
By (5.7)
\begin{lemma} At $p$, $0\leq D, A,E\leq \frac{4}{3}\|B\|^2$.
Furthermore, $D=0$
iff $A(\mbox{or}~E) =\frac{4}{3}\|B\|^2$  iff
$B=0$.  If $A=0$ ( $E=0$ resp.)  then $\wedge^2 B$ is selfdual 
(resp. antiselfdual).
\\[-2mm]
\end{lemma}
\noindent
Now we investigate when $\tilde{Q}_{\Omega}(B)\geq \delta\|B\|^2$.
Using the matrix (5.9)
\begin{eqnarray}
\lefteqn{\cos\theta \tilde{Q}_{\Omega}(B) = \cos\theta\|B\|^2
 -2\langle \Psi_+,\wedge^+_+B\rangle
-2\langle \Psi_-,\wedge^-_-B\rangle}\\[1mm]
&=&-\sm{\sum_{k}}\La{(}\sm{\sum_{r=1,2,3}}
~\sm{\frac{4}{3}}\langle \Lambda^{+}_r,
\wedge^2l_k(\Lambda^{+}_r)\rangle -
\|l_k\|^2\La{)}\nonumber\\[-1mm]
&&-s^2\sm{\sum_k}\La{(}\sm{\frac{4}{3}\langle \Lambda^{+}_1,
\wedge^2l_k(\Lambda^{+}_1)\rangle
-\frac{4}{3}\langle \Lambda^{+}_2,
\wedge^2l_k(\Lambda^{+}_2)\rangle
-\frac{4}{3}\langle \Lambda^{+}_3,
\wedge^2l_k(\Lambda^{+}_3)\rangle} \nonumber\\[-1mm]
&&~~~~~~~\sm{+\frac{4}{3}\langle \Lambda^{-}_1,
\wedge^2l_k(\Lambda^{-}_1)\rangle
+\frac{4}{3}\langle \Lambda^{-}_2,
\wedge^2l_k(\Lambda^{-}_2)\rangle
-\frac{4}{3}\langle \Lambda^{-}_3,
\wedge^2l_k(\Lambda^{-}_3)\rangle +
\frac{2}{3}\|l_k\|^2}\La{)}\nonumber \\[1mm]
&=&(D+ s^2E)-s^2(A+\sm{\frac{2}{3}}\|B\|^2).
\end{eqnarray}
\begin{lemma} Assume at $p$, for each $k$, $\|H^+l_k\|\leq \epsilon \|l_k\|$
where $0\leq \epsilon\leq 1$, and  $0\leq \tau 
\leq\frac{4(1-\epsilon)}{9}$ such that 
 $(1-\cos\theta)\leq \tau$.
 Then, at $p$, $\cos\theta\, \tilde{Q}_{\Omega}(B)\geq \delta \|B\|^2$ where 
$\delta=\frac{4(1-\epsilon)-9\tau}{3}\geq 0$.
\end{lemma}
\noindent
Note that we assume  $\cos\theta\geq \frac{5}{9}$.
Lemma 5.5 includes the case $\cos\theta\equiv 1$, that
implies $l_k$ hypercomplex, giving $\epsilon=1$
and $\delta=0$, as in proof of proposition 5.7. 
\\[2mm]
\em Proof. \em 
The condition on $H^+l_k$ implies 
$|\langle Q^+l_k,l_k\rangle +\frac{1}{3}|l_k|^2|\leq 
\frac{4}{3}\epsilon |l_k|^2$.
Then $D=\sum_k \frac{4}{3}\|l_k\|^2-( \frac{1}{3}\|l_k\|^2+
\langle Q^+l_k,l_k\rangle)\geq \frac{4}{3}(1- \epsilon)\|B\|^2$.
Using the bounds in lemma 5.4 and (5.15), we obtain $
\cos\theta\, \tilde{Q}_{\Omega}(B)\geq \frac{4(1- \epsilon)-{6}s^2}{3}
\|B\|^2$, that proves the lemma.\qed \\[3mm]
 Note that
(1.7) can be satisfied.
If $\oM$ is a quaternionic space form
of reduced scalar curvature $\nu={s^{\oM}}/{32}$, then
{
$$
\bar{R}(X,Y,Z,W)=\sm{\frac{\nu}{4}}\la{(}\langle
X\wedge Y, Z\wedge W\rangle +\sm{\sum}_r\langle J_rX\wedge J_rY, Z\wedge W\rangle
+\langle J_rX,Y\rangle\langle J_rZ,W\rangle\la{)}$$
}
and so
\begin{equation}
\sm{\sum}_{ki}~\bar{R}(X_i,X_k,X_i,\Phi(X_k))=9\nu (\-\cos\theta)(\cos\theta
-\sm{\frac{1}{3}})=4\nu s^2c^2
\end{equation}
\begin{proposition} Assume 
 $F:M\to \oM$ is a closed almost complex immersed 
submanifold such that (1.7) holds. \\[1mm]
(1) If  there exist constants  $0\leq \epsilon\leq 1$, 
$0\leq \tau <{4(1-\epsilon)}/{9}$, such that
 $(1-\cos\theta)\leq \tau$ and at each point $p\in M$ there exist
canonical frames $B$ and ${B'}^{\bot}$ such that
$\|H^+l_k\|\leq \epsilon \|l_k\|$, then $F$ is totally geodesic.
\\[1mm]
Furthermore, if $\oM$ is a quaternionic space form  then:\\[1mm]
(2) If $F:M\ra \oM$ is  parallel and $\nu\neq 0$ then either
$F$ is a quaternionic submanifold or a totally complex submanifold.\\[1mm]
(3) If $M$ is closed, $F$ has parallel mean curvature and $\wedge^2 B$
is anti-self-dual then $F$ is totally geodesic and if $\nu>0$
then either $F$ is totally complex or a quaternionic submanifold.
\end{proposition}
\noindent
\em Proof. \em 
(1) follows from previous lemma and theorem 1.5(A). If
we assume $\tau <\frac{4(1-\epsilon)}{9}$, it guarantees
$\delta>0$.
(2) From the proof of proposition 3.2, 
$\sum_k( \bar{R}(X_k,X_i)X_k)^{\bot}$ $=0$ and by (5.16),
$4\nu s^2c^2 =9\nu (1-\cos\theta)(\cos\theta -\frac{1}{3})=0$. (3)
We have $\Delta\cos\theta =-4\nu\, s^2c^2
-\cos\theta \|B\|^2$ $\leq 0$,
what implies $\cos\theta\geq  \frac{1}{3}$ is constant and so 
$-4\nu\, s^2c^2-\cos\theta\|B\|^2=0$.
\qed
\begin{proposition} 
Assume $M$ is a closed almost  complex submanifold with parallel mean
curvature  on a quaternionic space form $\oM$. Then $F$ is
totally geodesic if (1) or (2) below holds:\\[2mm]
(1)~ $\nu> 0$ and
$ \|B\|^2\leq 3\nu (\cos\theta-\frac{1}{3})$\\
(2)~ $\nu< 0$  and
$\|B\|^2\leq-\frac{27}{8}\nu (\cos\theta-\frac{1}{3})(
1-\cos\theta)$.
\end{proposition}
\noindent
\em Proof. \em 
We may write $\Delta \cos\theta$ given in (5.10) and using (5.15)  as
$$\Delta s^2 = 6\nu s^2c^2 + 4sc\sm{\sum}_j\langle \lnabo{X_j}H, U_j\rangle
-\sm{\frac{3}{2}}s^2(A+\sm{\frac{2}{3}}\|B\|^2)+\sm{\frac{3}{2}}(D+ s^2E).$$
(1) The conditions  imply
$\frac{3}{2}s^2(A+\frac{2}{3}\|B\|^2)
\leq 6\nu s^2c^2$ and so $\Delta s^2\geq 0$.
 (2)
Under the assumptions, 
 $\frac{3}{2}(D+s^2E)\leq 2(1+ s^2)\|B\|^2\leq 4\|B\|^2\leq -6s^2c^2\nu$, and
we have $\Delta s^2 \leq 0$. In both cases (1)(2)
 we conclude that $s$ is constant, and again
that $\Delta s^2=0$. 
This implies  in case (1) $D+s^2E=0$, and in case (2) 
$s^2(A+ \frac{2}{3}\|B\|^2)=0$, what leads to the conclusion in
the proposition (see lemma 5.4).\qed
\subsection{The special calibrations}
{\em The special Lagrangian calibration.}~ 
Let $(\oM,g,J,\rho)$ be a Calabi-Yau manifold of complex dimension $k$
with  holomorphic volume element $\rho\in\wedge^{(k,0)}M$.
Then $Re(\rho)$ is the Lagrangian calibration
and calibrates the special Lagrangian submanifolds.
On $\oM$ it is also defined the K\"{a}hler calibrations. If $k=4$,
there is also
a $S^1$-family of Cayley calibrations
$\Omega_{\theta}=-\frac{1}{2}w^2+Re(e^{\theta}\rho)$, that calibrates
the Cayley $4$-submanifolds.\\[2mm]
{\em The Cayley calibration.}~
If $(\oM^8,\bar{g}, \Omega)$ is a $Spin(7)$ 8-dimensional manifold, 
then it is defined
a Cayley calibration $\Omega$.  Given a  $Spin(7)$-frame $e_i$ that identifies
$T_p\oM$ with the space of octonions $\R^8$,  $\Omega$ is 
the 4-form defined by
$\Omega(x,y,z,w)=\langle x,y\times z\times w \rangle $, where
the cross product of three vectors is defined in $\R^8$.
A Calabi-Yau 4-fold is also a $Spin(7)$ manifold and any Cayley
calibrations defined above corresponds to this definition. 
If $F:M^4\ra \oM$ is an immersed 4-submanifold, then
$\Phi(X_1)=(X_2 \times X_3\times X_4)^{\bot}$, where $X_i$ is
a d.o.n. basis of $T_pM$.\\[1mm]
\noindent
{\em The associative and the co-associative calibration.}~
Let $(\oM^7,\bar{g},\phi)$ be a $G_2$ Riemannian manifold with
a closed $G_2$ 3-form $\phi$. Identifying
$T_p\oM$ with $\mathbb{R}^7=Im (\mathbb{R}^8)$ by a $G_2$-frame, 
$\phi(x,y,z)=\langle x,yz\rangle$ where on the right hand side
it is considered the octonion product. This is the associative
calibration. The co-associative calibration is $\psi=*\phi$
and satisfies $\psi(x,y,z,w)=\frac{1}{2}\langle x,[y,z,w]\rangle $
where $[y,z,w]=(yz)w-y(zw)$ is the associator operator. The forms
$\phi$ and $\psi$  calibrate respectively the
associative 3-dimensional submanifolds and the co-associative 4-dimensional 
submanifolds. If $F:M\ra \oM$ is an immersed 3-submanifold, 
$\Phi_{\phi}(X_1)=(X_2X_3)^{\bot}$, where $X_1,X_2,X_3$ is any d.o.n. basis of 
$T_pM$. If $F$ is an immersed 4-submanifold,
$\Phi_{\psi}(X_1)=[X_2,X_3,X_4]^{\bot}$.
If $N$ is a Calabi-Yau 3-fold, then $N\times S^1$ or $N\times \R$
are $G_2$-manifolds with $\phi= 1^*\wedge w+Re(\rho)$
and $\psi=\frac{1}{2}w\wedge w-1^*\wedge Im(\rho)$.
If $N$ is a $G_2$ manifold, then $N\times S^1$ or $N\times \R$ with
 $\Omega=1^*\wedge \phi +\psi$   are
$Spin(7)$ manifolds.

\section{Acknowledgements}
The authors would like to thank the  referee for helpful
comments and valuable suggestions.
The second author would like to thank Dmitri Alekseevsky for
many conversations in quaternionic K\"{a}hler geometry.
The first author is partially supported by NSFC (No.10501011) and by
Funda\c{c}\~{a}o Ci\^{e}ncia e Tecnologia (FCT) through a FCT
fellowship SFRH/BPD/26554/2006. The second author is
partially supported by FCT through the Plurianual of CFIF
and POCI-PPCDT/MAT/60671/2004.\\[6mm]
\noindent
(1)~\begin{minipage}[t]{13cm} School of Mathematics and Computer Science,\\
Hubei University, Wuhan, 430062,
P. R. China;\\
 liguanghan@163.com
\end{minipage}
\\[3mm]
(2)~\begin{minipage}[t]{15cm}
Centro de F\'{\i}sica das Interac\c{c}\~{o}es Fundamentais, \\
Instituto Superior T\'{e}cnico,
 Technical University of Lisbon, \\
Edif\'{\i}cio Ci\^{e}ncia, Piso 3,
Av.\ Rovisco Pais,
1049-001 Lisboa, Portugal;~\\
 isabel.salavessa@ist.utl.pt
\end{minipage}

\end{document}